\documentclass[12pt,french]{article}

\usepackage{amsthm,amsmath,amssymb,amsbsy,bbm,mathrsfs,supertabular,eurosym,
	graphicx,enumitem,xcolor}
\usepackage[authoryear]{natbib}
\usepackage[french]{babel}
\usepackage{ae}			
\usepackage[latin1]{inputenc} 	
\usepackage[T1]{fontenc}

\newtheorem{thm}{Th\'eor\`eme}

\newtheorem{prop}{Proposition}

\newtheorem{df}{D\'efinition}

\newtheorem{rem}{Remarque\!}

\newcommand{\pa}[1]{\left({#1}\right)}
\newcommand{\norm}[1]{\left\|{#1}\right\|}
\newcommand{\cro}[1]{\left[{#1}\right]}
\newcommand{\ab}[1]{\left|{#1}\right|}
\newcommand{\ac}[1]{\left\{{#1}\right\}}

\def\argmin{\mathop{\rm argmin}}
\def\argmax{\mathop{\rm argmax}}

\def\pen{\mathop{\rm pen}\nolimits}
\def\Var{\mathop{\rm Var}\nolimits}

\newcommand{\fleche}{\mathop{\longrightarrow}}
\newcommand{\dfleche}[1]{\,\displaystyle{\mathop{\longrightarrow}_{#1}}\,}

\newcommand{\CV}[1]{\dfleche{#1}}

\def\E{{\mathbb{E}}}

\def\N{{\mathbb{N}}}
\def\P{{\mathbb{P}}}
\def\Q{{\mathbb{Q}}} 
\def\R{{\mathbb{R}}}

\def\sF{{\mathscr{F}}}

\def\sP{{\mathscr{P}}}
\def\sQ{{\mathscr{Q}}}

\def\sW{{\mathscr{W}}}
\def\sX{{\mathscr{X}}}

\DeclareMathAlphabet{\mathscrbf}{OMS}{mdugm}{b}{n}

\def\sbF{{\mathscrbf{F}}}

\def\sbP{{\mathscrbf{P}}}
\def\sbQ{{\mathscrbf{Q}}} 

\def\sbS{{\mathscrbf{S}}}

\def\sbX{{\mathscrbf{X}}}

\def\cF{{\mathcal{F}}}

\def\cK{{\mathcal{K}}}
 
\def\cM{{\mathcal{M}}}
\def\cN{{\mathcal{N}}}

\def\cP{{\mathcal{P}}}

\def\cU{{\mathcal{U}}}

\def\gp{{\mathbf{p}}}
\def\gq{{\mathbf{q}}}

\def\gP{{\mathbf{P}}}
\def\gQ{{\mathbf{Q}}} 
\def\gR{{\mathbf{R}}}
 
\def\gT{{\mathbf{T}}}

\def\gX{{\mathbf{X}}}

\newcommand{\bs}[1]{\boldsymbol{#1}}

\def\bsX{{\bs{X}}}
\def\biX{{\bs{X}}}

\def\gmu{{\bs{\mu}}}

\def\gup{\bs{\Upsilon}}

\newlist{lista}{enumerate}{1}
\setlist[lista,1]{label=\alph*),ref=\alph*)}

\newlist{listi}{enumerate}{1}
\setlist[listi,1]{label=(\roman*),ref=(\roman*),align=left}

\newcommand{\eref}[1]{(\ref{#1})}
\def\1{1\hskip-2.6pt{\rm l}}
\def\<{{\langle}}
\def\>{{\rangle}}

\newcommand{\etc}[1]{#1_1,\ldots,#1_n}

\newcommand{\st}{\strut}

\newcommand{\ord}[1]{$^{\mbox{#1}}$}

\def\et{^{\star}}
\def\eps{{\varepsilon}}

\newcommand{\Vx}[1]{V_{#1}}   

\begin{document}
\title{\mbox{}\vspace{-10mm}\\
{\Huge Une alternative robuste au maximum de vraisemblance: la $\rho$-estimation}
}
\author{{\LARGE{\bf Yannick Baraud et Lucien Birg\'e}}}

\date{1\ord{er} juillet 2017}
\maketitle

\begin{abstract}
Cet article est fond\'e sur les notes du mini-cours que nous avons donn\'e le 5 janvier 2017 \`a l'Institut Henri Poincar\'e \`a l'occasion d'une journ\'ee organis\'ee par la Soci\'et\'e Fran\c caise de Statistique et consacr\'ee \`a la Statistique Math\'ematique. Il vise \`a donner un aper\c{c}u de la m\'ethode de $\rho$-estimation comme alternative \`a celle du maximum de vraisemblance, ainsi que des propri\'et\'es d'optimalit\'e et de robustesse des $\rho$-estimateurs. Cette m\'ethode s'inscrit dans une longue lign\'ee de recherche initi\'ee par de c\'el\`ebres statisticiens tels que Sir Ronald Fisher, avec le maximum de vraisemblance dans les ann\'ees 20, ou Lucien Le Cam, avec des estimateurs fond\'es sur des tests entre boules de Hellinger dans les ann\'ees 70, et dont l'objectif a \'et\'e de produire des m\'ethodes d'estimation poss\'edant de bonnes propri\'et\'es pour un ensemble de cadres statistiques aussi vaste que possible. Plus r\'ecemment, Birg\'e avec les d- puis T-estimateurs, a \'etendu les r\'esultats de Le Cam dans diverses directions, en particulier la robustesse et l'adaptation. Nous mettrons en lumi\`ere les liens forts qui existent entre les $\rho$-estimateurs et ces pr\'ed\'ecesseurs, notamment les estimateurs du maximum de vraisemblance, mais montrerons \'egalement, au travers d'exemples choisis, que les $\rho$-estimateurs les surpassent par bien des aspects. 
\end{abstract}

\section{Bref historique\label{I}}

Il y a bien longtemps que les statisticiens ont cherch\'e une m\'ethode g\'en\'erale pour construire, dans des contextes divers, des estimateurs disposant de bonnes propri\'et\'es, en particulier asymptotiques. En r\'egression, il y a eu le minimum $\ell_{1}$ attribu\'e \`a Pierre Simon de Laplace et les moindres carr\'es (minimum $\ell_{2}$) de Carl Friedrich Gau\ss\ autour des ann\'ees 1800. Pour des mod\`eles param\'etriques g\'en\'eraux sur variables i.i.d.\ la m\'ethode des moments est connue depuis longtemps et celle du maximum de vraisemblance, d\'evelopp\'ee et popularis\'ee d'abord par Sir Ronald Fisher dans les ann\'ees 1920, a connu un succ\`es consid\'erable et elle est, aujourd'hui encore, tr\`es largement utilis\'ee. L'estimateur du maximum de vraisemblance (en abr\'eg\'e e.m.v.) est, en un certain sens, le premier estimateur \`a vis\'ees ``universelles" dans la mesure o\`u il s'applique \`a de tr\`es nombreux mod\`eles statistiques et permet aussi de retrouver, en r\'egression, la m\'ethode des moindres carr\'es (maximum de vraisemblance pour des erreurs gaussiennes) ou du minimum $\ell_{1}$ (lorsque les erreurs ont une loi de Laplace).

L'e.m.v., qui a d'excellentes propri\'et\'es asymptotiques dans les mod\`eles param\'etriques r\'eguliers, sous r\'eserve de consistance, souffre n\'eanmoins, dans un cadre g\'en\'eral, de certains d\'efauts connus depuis longtemps et mis en \'evidence par des statisticiens tels que Bahadur d\`es 1958 --- cf. Bahadur~\citeyearpar{Bahadur58} --- ou Le~Cam --- cf.\ Le~Cam~\citeyearpar{Lecam-MLE} ---, entre autres. Il n'est pas possible de construire une th\'eorie g\'en\'erale de consistance des estimateurs du maximum de vraisemblance sans utiliser des hypoth\`eses relativement fortes comme on peut le voir par exemple dans le livre de van~der~Vaart~\citeyearpar{MR1652247}. Le probl\`eme est encore plus d\'elicat si l'on veut \'etudier les vitesses de convergence des e.m.v. Le Cam en \'etait tr\`es conscient, cela appara\^it dans ses \'ecrits et il s'est efforc\'e de r\'esoudre le probl\`eme en construisant --- cf.\ Le~Cam~\citeyearpar{MR0334381} et \citeyearpar{MR0395005} --- un estimateur qui est automatiquement consistant sous des hypoth\`eses relativement faibles. C'\'etait l\`a une nouvelle tentative, plus f\'econde, pour fabriquer une m\'ethode d'estimation ``universelle", au moins pour les mod\`eles en variables ind\'ependantes.

C'est en se fondant sur ces travaux de Le Cam en vue de g\'en\'eraliser les choses \`a des variables non n\'ecessairement ind\'ependantes, en recherchant des vitesses de convergence et en \'etudiant leur optimalit\'e \'eventuelle et leur robustesse, que le second auteur de cet article a construit ce qu'il a appel\'e les d-estimateurs dans Birg{\'e}~\citeyearpar{MR722129}. Plus tard ceux-ci se sont transform\'es en T-estimateurs par l'adjonction de p\'enalit\'es qui ont conduit \`a la s\'election de mod\`ele et l'adaptativit\'e dans Birg{\'e}~\citeyearpar{MR2219712}. N\'eanmoins, les T-estimateurs, tout comme leurs pr\'ed\'ecesseurs, souffrent encore de diverses limitations, en particulier la n\'ecessit\'e de travailler avec des mod\`eles compacts, ce qui s'av\`ere bien trop restrictif pour les applications \`a la r\'egression.

Il est important de noter ici que ces diff\'erentes constructions, celle de Le Cam comme celles  de Birg\'e, reposent sur la construction de familles de tests robustes (tests entre deux boules de Hellinger initialement). On peut d'ailleurs aussi interpr\'eter le maximum de vraisemblance comme le r\'esultat de la combinaison de nombreux tests de rapport de vraisemblances comme nous le verrons plus loin. C'est en partant d'une id\'ee originale du premier auteur pour construire de nouvelles familles de tests --- cf. Baraud~\citeyearpar{MR2834722} --- qu'est venue l'id\'ee des $\rho$-estimateurs. Ces nouveaux tests sont fond\'es sur l'id\'ee suivante : \'etant donn\'e trois lois de probabilit\'e $P$, $Q$ et $R$ et $\etc{X}$ i.i.d.\ de loi $P$, on peut d\'eduire un test entre $Q$ et $R$ (comme meilleur approximant de $P$) d'une estimation du signe de la diff\'erence
$h^{2}(P,Q)-h^{2}(P,R)$ o\`u $h$ d\'esigne la distance de Hellinger.

Dans le cas de mod\`eles compacts, les $\rho$-estimateurs peuvent \^etre vus comme une variante des T-estimateurs fond\'es \'egalement sur des tests robustes, mais ils permettent en outre de traiter des probl\`emes \`a mod\`eles non compacts comme la r\'egression lin\'eaire en plan fixe ou al\'eatoire, avec des erreurs de lois tr\`es vari\'ees, ou des probl\`emes d'estimation de densit\'es sous contraintes de forme. Comme les T-estimateurs ils sont robustes et permettent de consid\'erer de nombreux mod\`eles simultan\'ement, d'o\`u leurs propri\'et\'es d'adaptation, mais leurs performances vont tr\`es au-del\`a et permettent bien davantage que les T-estimateurs, par exemple l'estimation d'une densit\'e unimodale sur $\R$ de mode et support inconnus. Enfin, comme nous le verrons, dans des situations favorables, l'e.m.v.\ est un $\rho$-estimateur ; la $\rho$-estimation peut donc, dans une certaine mesure, \^etre consid\'er\'ee comme une mani\`ere de rendre robuste la c\'el\`ebre  proc\'edure introduite il y a pr\`es d'un si\`ecle par Sir Ronald Fisher.

Ce qui va suivre peut \^etre consid\'er\'e comme une introduction (voire une incitation) \`a la lecture des articles suivants :
Baraud, Birg\'e et Sart~\citeyearpar{MR3595933}, Baraud et Birg\'e~(2016a) et Baraud et Birg\'e~(2016b).\\
\phantom{\citeyearpar{MR3565484}} 
\phantom{\citeyearpar{Baraud:2016kq}}

\section{Quelques probl\`emes li\'es \`a l'utilisation des estimateurs du maximum de vraisemblance \label{EMV}}

\subsection{Instabilit\'e de la m\'ethode du maximum de vraisemblance\label{EMV1}}
On dispose d'un million d'observations $X_{1},\ldots,X_{10^{6}}$ que l'on suppose  i.i.d.\ de loi uniforme sur $[0,\theta]$, de param\`etre inconnu $\theta>0$ et que l'on mod\'elise donc comme telles. L'e.m.v.\ relatif \`a ce mod\`ele est $\widehat \theta=\max_{i=1,\ldots,10^{6}}X_{i}$. Si toutes les observations suivent effectivement ce mod\`ele avec $\theta=1$ sauf une seule d'entre elles qui, suite \`a une erreur de manipulation ou de transcription,  vaut 100, alors le maximum de vraisemblance vaut 100 et l'erreur d'estimation de $\theta$ sera d'au moins 99. Il suffit donc d'une seule observation aberrante parmi un million pour rendre l'e.m.v.\ inop\'erant dans ce cadre. 

Le m\^eme probl\`eme se pose pour l'estimation du param\`etre de centrage $\theta$ de la gaussienne $\cN(\theta,1)$ par la moyenne empirique, qui est aussi l'e.m.v.\ ; une seule observation aberrante peut conduire \`a une erreur d'estimation consid\'erable, autrement dit l'e.m.v.\ n'est absolument pas {\em robuste}. Pour plus de d\'etails sur la notion de {\em robustesse}, on pourra se r\'ef\'erer \`a l'excellent livre de Peter Huber --- Huber~\citeyearpar{MR606374} --- et \`a sa bibliographie.

La d\'etection d'une valeur aberrante en dimension 1 (ou 2) est en g\'en\'eral relativement facile car on peut visualiser les donn\'ees. En revanche, les choses se compliquent s\'erieusement en grande dimension !

\subsection{Non-existence de l'estimateur du maximum de vraisemblance\label{EMV2}}
Soit la densit\'e $p$ sur $\R$ d\'efinie par 
\begin{equation}
p(x)=(1/6)\left[|x|^{-1/2}\1_{(0,1]}(|x|)+x^{-2}\1_{]1;+\infty[}(|x|)\right],
\label{eq-densitp}
\end{equation}
laquelle est sym\'etrique, mais discontinue et non born\'ee en 0. Si l'on consid\`ere la
famille de translation $\{p_\theta(x), \theta\in\Theta\}$ o\`u $p_\theta(x)=p(x-\theta)$ et $\Theta$ est un intervalle compact de $\R$, la vraisemblance s'\'ecrit $p(X_{1}-\theta)\times\ldots\times p(X_{n}-\theta)$ et tend vers l'infini d\`es que $\theta\rightarrow X_i$, quel que soit $i$. On ne peut
donc d\'efinir un e.m.v. Le m\^eme ph\'enom\`ene se produit pour tout mod\`ele de translation pour lequel la densit\'e $p$ n'est pas born\'ee.

Dans le cas o\`u $p$ est donn\'ee par \eref{eq-densitp}, il est n\'eanmoins facile de montrer que la m\'ediane empirique va converger vers $\theta$ \`a vitesse $1/n$. On peut m\^eme construire des intervalles de confiance non-asymptotiques pour $\theta$; ce n'est donc pas le probl\`eme qui est difficile \`a r\'esoudre mais plut\^ot  la m\'ethode du maximum de vraisemblance qui est inadapt\'ee \`a celui-ci.

N\'eanmoins, si l'on remplace le mod\`ele initial $\{p_\theta, \theta\in\Theta\}$ par le mod\`ele approch\'e $\{p_\theta, \theta\in\Theta'\}$ o\`u $\Theta'$ est une discr\'etisation finie de l'intervalle $\Theta$, l'e.m.v.\ relatif \`a $\Theta'$ sera bien d\'efini.

\subsection{Le risque quadratique de l'e.m.v.\ peut \^etre bien plus grand que le risque minimax}\label{sect-EMVbad}
Soit $\biX= (X_0,\ldots,X_k)$ un vecteur gaussien $(k+1)$-dimensionnel de loi ${\cal N}(\theta,I_{k+1})$, o\`u $I_{k+1}$ est la matrice identit\'e de dimension $k+1$. \'Etant donn\'e un vecteur $\theta=(\theta_0,\ldots,\theta_k)$ de $\Bbb{R}^{k+1}$, on note $\theta'=(\theta_1,\ldots,\theta_k)$ sa projection sur le sous-espace vectoriel de dimension $k$ engendr\'e par les $k$ derni\`eres coordonn\'ees, $\|\theta\|$ sa norme euclidienne et $\E_{\theta}$ l'esp\'erance quand $\theta$ est la vraie valeur du param\`etre.
%
\begin{prop}\label{GM}
Supposons que $k\ge128$ et consid\'erons l'espace de param\`etres
\[
\Theta=\left\{\left.\theta\in\Bbb{R}^{k+1}\,\st\right|\,|\theta_0|\le k^{1/4}
\quad\mbox{et}\quad\|\theta'\|\le2\left(1-k^{-1/4}|\theta_0|\right)\right\}.
\]
Le risque quadratique de l'e.m.v.\ $\widehat{\theta}$ et le risque
minimax sur $\Theta$ v\'erifient respectivement
\begin{equation}
\sup_{\theta\in\Theta}\Bbb{E}_\theta\left[\|\theta-\widehat{\theta}\|^2\right]\ge(3/4)\sqrt{k}+3\qquad\mbox{et}\qquad\inf_{\widetilde{\theta}}\sup_{\theta\in\Theta}
\Bbb{E}_\theta\left[\|\theta-\widetilde{\theta}\|^2\right]\le5.
\label{eq-RMLE-Gau}
\end{equation}
\end{prop}
On voit que le risque quadratique de l'e.m.v.\ sur $\Theta$ peut \^etre beaucoup plus grand que le risque minimax. La d\'emonstration se trouve dans Birg{\'e}~\citeyearpar{MR2219712}.

Pour obtenir la borne sup\'erieure du risque quadratique dans \eref{eq-RMLE-Gau} il suffit de remplacer le mod\`ele initial $\Theta$ par le mod\`ele approch\'e $\Theta'=\ac{\left.(\theta_{0},0)\in\R^{k+1}\,\st\right|\, \theta_{0}\in\R}$ et de construire l'e.m.v.\ sur $\Theta'$ ce qui donne $\widetilde \theta=(\widetilde \theta_{0},0)=(X_{0},0)$. On v\'erifie alors que
\[
\Bbb{E}_\theta\left[\|\theta-\widetilde{\theta}\|^2\right]=\E_{0}\cro{(\theta_{0}-X_{0})^{2}}+\|\theta'\|^2\le 1+4=5\quad\mbox{pour tout }\theta\in\Theta.
\]

\subsection{La valeur et les performances de l'e.m.v.\ d\'ependent du choix de la famille des densit\'es}
Que la valeur de l'e.m.v.\ d\'epende du choix que l'on fait des densit\'es par rapport \`a une mesure dominante est \`a peu pr\`es \'evident et connu depuis longtemps. Qu'un choix inadapt\'e puisse mener \`a des r\'esultats catastrophiques, comme le montre l'exemple suivant, l'est sans doute nettement moins.
%
\begin{prop}\label{prop-P1}
Soit une suite $(X_{k})_{k\ge 1}$ de v.a.r.\ i.i.d.\ de loi $P_{\theta}=\cN(\theta,1)$ pour un r\'eel $\theta$ inconnu. Pour tout $\theta\in\R$, fixons la version $p_{\theta}$ de la densit\'e $dP_{\theta}/dP_{0}$ de $P_{\theta}$ par rapport \`a la mesure de r\'ef\'erence $P_{0}=\cN(0,1)$ de la mani\`ere suivante :
%
\begin{equation}
p_{\theta}(x)=
\left\{
\begin{array}{ll}
\exp\cro{\theta x-\left(\theta^{2}/2\right)} &\mbox{si $x\ne \theta$ ou $\theta\le 0$},\\
\exp\cro{\theta x-\left(\theta^{2}/2\right)+(\theta^{2}/2)\exp\left(x^{2}\right)}\; & \mbox{sinon.}
\end{array}
\right.
\label{eq-bad MLE}
\end{equation}
Alors, quelle que soit la vraie valeur de $\theta\in\R$, sur un ensemble dont la $P_{\theta}$-probabilit\'e tend vers 1 quand $n$ tend vers l'infini, l'e.m.v.\ s'\'ecrit $\widehat{\theta}_{n}=X_{(n)}=\max\{X_{1},\ldots,X_{n}\}$ et n'est donc pas consistant.
\end{prop}
On trouvera la preuve de cette assertion dans Baraud et Birg\'e~(2016b).

Si l'on remplace le mod\`ele initial avec $\theta\in\R$ par sa restriction \`a $\theta\in\Q$ et si l'on accepte de prendre pour maximum de vraisemblance un point maximisant la vraisemblance appproximativement seulement, le probl\`eme ne se posera plus.

\subsection{Construction de l'e.m.v.\ sur un mod\`ele approch\'e}\label{ex-sect2.5}
Nous venous de remarquer que certains des probl\`emes rencontr\'es par l'e.m.v.\  pouvaient \^etre r\'esolus en construisant l'estimateur non pas sur l'ensemble initial des param\`etres $\Theta$ mais sur un sous-ensemble $\Theta'$ de $\Theta$ ayant de bonnes qualit\'es d'approximation par rapport aux \'elements de $\Theta$. Maleureusement, l'utilisation de tels sous-ensembles approchants peut aussi conduire \`a d'autres ennuis comme le montre l'exemple suivant.

Consid\'erons le probl\`eme d'estimation du param\`etre $\theta\in\R$ \`a partir de l'observation d'un $n$-\'echantillon $X_{1},\ldots,X_{n}$ dont la loi appartient \`a l'ensemble 
de probabilit\'es $\sP_{\eps}=\ac{P_{\alpha,\theta},\; (\alpha,\theta)\in [0,\eps]\times\R}$ avec $0\le\varepsilon<1/2$ et
\[
P_{\alpha,\theta}=(1-\alpha)\,\cU([\theta,\theta+1])+\alpha\,\cU([100+\theta,101+\theta],
\]
o\`u $\cU([a,b])$ d\'esigne la loi uniforme sur l'intervalle $[a,b]$ pour $a<b$. Si $\varepsilon$ est petit le mod\`ele plus simple $\sP_{0}\subset \sP_{\eps}$ des lois uniformes sur $[\theta,\theta+1]$ avec $\theta\in\R$ fournit une bonne approximation de $\sP_{\eps}$. Comme la loi $P_{\alpha,\theta}$ est un m\'elange de deux lois uniformes sur deux intervalles disjoints, la vraisemblance sur le mod\`ele $\sP_{0}$ est constamment nulle sauf si toutes les observations appartiennent \`a l'un ou l'autre de ces deux intervalles ce qui arrive avec une probabilit\'e 
\[
(1-\alpha)^{n}+\alpha^{n}\le\exp[-n\alpha]+\alpha^{n}.
\]
Lorsque $\alpha$ est grand devant $1/n$ et petit devant $1/2$, cette probabilit\'e est petite. Donc, en dehors de cet \'ev\'enement de petite probabilit\'e, la m\'ethode du maximum de vraisemblance sur le mod\`ele approch\'e $\sP_{0}$ ne fournit aucune estimation de $\theta$.

\subsection{Les mod\`eles approch\'es\label{S2.6}}
Si nous consid\'erons les exemples pr\'ec\'edents, nous pouvons constater que les probl\`emes li\'es \`a l'utilisation de l'e.m.v. sur un mod\`ele statistique donn\'e par une famille de lois $\{P_{\theta},\theta\in \Theta\}$ peuvent souvent \^etre en partie r\'esolus en rempla\c cant le mod\`ele initial par un mod\`ele approch\'e $\{P_{\theta},\theta\in\Theta'\}$. Pour obtenir un r\'esultat g\'en\'eral fond\'e sur l'estimation sur des mod\`eles approch\'es, il convient de s'assurer que, si l'on utilise comme mod\`ele une famille de lois $\cP=\{P_{\theta},\theta\in \Theta'\}$ et si la vraie loi $P$ des observations est proche d'un \'el\'ement de $\cP$, le comportement de l'estimateur construit sur $\cP$ est \`a peu pr\`es le m\^eme que ce qu'il serait si $P$ appartenait effectivement \`a $\cP$, c'est-\`a-dire que l'estimateur est {\em robuste}. Malheureusement l'e.m.v.\ ne l'est pas toujours, d'o\`u la n\'ecessit\'e de le remplacer par un estimateur robuste permettant de travailler avec des mod\`eles approch\'es de mani\`ere syst\'ematique.

\section{Comment r\'esoudre les probl\`emes pr\'ec\'edents : une heuristique\label{S3}}
Pla\c cons-nous ici dans le cadre d'un $n$-\'echantillon $\biX=(\etc{X})$ et d'une famille $\sbS$ de densit\'es par rapport \`a une mesure de r\'ef\'erence $\mu$ avec $P_{t}=t\cdot\mu$ pour $t\in\sbS$. Notons $\P_{t}[A]$ la probabilit\'e de l'\'ev\`enement $A$ quand les $X_{i}$ ont la loi $P_{t}$, la vraisemblance de $t$ par 
\[
\Vx{\biX}(t)=\prod_{i=1}^{n}t(X_{i}),\quad\mbox{donc}\quad\log\left(\st\Vx{\biX}(t)\right)=\sum_{i=1}^{n}\log\left(t(X_{i})\st\right),
\]
et supposons, pour simplifier les choses, que $\P_{s}\left[\Vx{\biX}(t)=\Vx{\biX}(u)\right]=0$, quels que soient $s,t,u\in\sbS$ et $t\ne u$. Alors on sait que, si $P_{s}\ne P_{t}$
\begin{equation}
\P_{s}\left[\Vx{\biX}(t)>\Vx{\biX}(s)\right]\le\rho^{n}(P_{s},P_{t}),
\label{eq-test}
\end{equation}
avec
\[
\rho(P_{s},P_{t})=\int\sqrt{s(x)t(x)}\,d\mu(x)<1.
\]
La quantit\'e $\rho(P_{s},P_{t})$, d\'enomm\'ee {\em affinit\'e de Hellinger entre $P_{s}$ et $P_{t}$}, est seulement fonction des deux lois $P_{s}$ et $P_{t}$ et non pas du choix de la mesure dominante $\mu$ ni de celui des versions des densit\'es $dP_{s}/d\mu$ et $dP_{t}/d\mu$. Comme $\rho(P_{s},P_{t})<1$, $\rho^{n}(P_{s},P_{t})$ tend vers $0$ quand $n$ tend vers l'infini et l'on voit bien qu'en choisissant en fonction de $\biX$ celle des deux lois qui a la plus grande vraisemblance, on trouvera asymptotiquement la vraie loi des $X_{i}$. On contr\^ole m\^eme les probabilit\'es d'erreurs du test pour $n$ fix\'e :
\[
\max\left\{\P_{s}\left[\Vx{\biX}(t)>\Vx{\biX}(s)\right],\P_{t}\left[\Vx{\biX}(t)<\Vx{\biX}(s)\right]\right\}
\le\rho^{n}(P_{s},P_{t}),
\]
puisque $\Vx{\biX}(t)\ne\Vx{\biX}(s)$ presque s\^urement. Le r\'esultat s'\'etend imm\'ediatement au cas d'un mod\`ele $\sbS$ fini et identifiable. En effet, pour tout $s\in\sbS$,
\[
\P_{s}\left[\sup_{t\in\sbS}\Vx{\biX}(t)>\Vx{\biX}(s)\right]\le\sum_{t\in\sbS,\,t\ne s}\P_{s}\left[\Vx{\biX}(t)>\Vx{\biX}(s)\right]
\le\sum_{t\in\sbS,\,t\ne s}\rho^{n}(P_{s},P_{t}),
\]
donc, si
\[
\widehat{t}_{n}(\biX)=\argmax_{t\in\sbS}\Vx{\biX}(t)=\argmax_{t\in\sbS}
\left[\sum_{i=1}^{n}\log\left(t(X_{i})\st\right)\right],
\]
alors
\[
\sup_{s\in\sbS}\P_{s}\left[\widehat{t}_{n}(\biX)\ne s\right]\CV{n\rightarrow+\infty}0.
\]
Remarquons ici que la recherche de $\widehat{t}_{n}$ revient \`a faire simultan\'ement tous les tests de rapport de vraisemblances entre les paires $(t,u)$ de $\sbS^{2}$ avec $t\ne u$ en comparant $\Vx{\biX}(t)$ \`a $\Vx{\biX}(u)$, ce qui revient \`a \'etudier le signe de 
\[
\log\left(\frac{\Vx{\biX}(t)}{\Vx{\biX}(u)}\right)=\sum_{i=1}^{n}\log\left(\frac{t(X_{i})}{u(X_{i})}\right)
\]
et \`a choisir l'\'el\'ement de $\sbS$ accept\'e par tous ces tests simultan\'ement, lequel existe toujours sous nos hypoth\`eses.
 
Malheureusement, ce r\'esultat ne s'\'etend pas sans hypoth\`eses fortes au cas d'un ensemble $\sbS$ infini ni \`a celui d'une vraie loi $P\et$ des $X_{i}$ n'appartenant pas \`a $\{P_{t}, t\in\sbS\}$, m\^eme lorsque $P\et$ est tr\`es proche de l'ensemble en question. Une des raisons de cet \'echec est li\'ee au fait que la fonction $\log$ n'est pas born\'ee comme le montre l'exemple de la section~\ref{EMV2}. Les tests de rapport de vraisemblances \'evoqu\'es ci-dessus reposent sur l'\'evaluation des rapports $\Vx{\biX}(t)/\Vx{\biX}(u)$ ou, de mani\`ere \'equivalente, sur les signes des quantit\'es $\sum_{i=1}^{n}\log\left(t(X_{i})/u(X_{i})\st\right)$. Mais il faut noter que c'est la propri\'et\'e magique de la fonction $\log$ : $\log(a/b)=\log a-\log b$ qui fait que trouver un point $\widehat{t}_{n}$ de $\sbS$ tel que $\log\left[\Vx{\biX}(\widehat{t}_{n})/\Vx{\biX}(u)\right]>0$ pour tout $u\ne\widehat{t}_{n}$ \'equivaut \`a chercher le point $\widehat{t}_{n}$ qui maximise sur $\sbS$ la fonction $u\mapsto\sum_{i=1}^{n}\log\left(u(X_{i})\st\right)$. Les fonctions de la forme $\alpha\log$ avec $\alpha>0$ sont les seules qui poss\`edent cette propri\'et\'e mais elles ne sont pas born\'ees ce qui explique certains des probl\`emes que l'on rencontre avec la m\'ethode du maximum de vraisemblance.

N\'eanmoins, l'id\'ee m\^eme de rapport de vraisemblances nous dit que $P_{t}$ est d'autant plus \og vraisemblable \fg{} compar\'e \`a $P_{u}$ que le rapport $\Vx{\biX}(t)/\Vx{\biX}(u)$ est plus grand ou que $\sum_{i=1}^{n}\log\left(t(X_{i})/u(X_{i})\st\right)$ l'est. De mani\`ere heuristique, la valeur de 
\begin{equation}
\sum_{i=1}^{n}\log\left(\frac{t(X_{i})}{u(X_{i})}\right)=2\sum_{i=1}^{n}\log\left(\sqrt{\frac{t(X_{i})}{u(X_{i})}}\right)
\label{eq-logvrais}
\end{equation}
permet de comparer la qualit\'e de $t$ par rapport \`a celle de $u$. Soit alors $\psi$ une fonction strictement croissante, born\'ee et telle que $\psi(1)=0$. On peut envisager de comparer, de mani\`ere analogue, la qualit\'e de $t$ par rapport \`a celle de $u$ en rempla\c cant dans (\ref{eq-logvrais}),
\[
\sum_{i=1}^{n}\log\left(\sqrt{\frac{t(X_{i})}{u(X_{i})}}\right)\quad\mbox{par}
\quad\sum_{i=1}^{n}\psi\left(\sqrt{\frac{t(X_{i})}{u(X_{i})}}\right),
\]
de sorte que l'on puisse dire que la loi $P_{u}$ est d'autant moins bonne au vu des observations $\etc{X}$ que
\[
\sup_{t\in\sbS}\sum_{i=1}^{n}\psi\left(\sqrt{\frac{t(X_{i})}{u(X_{i})}}\right)\,\mbox{ est grand}.
\]
Comme le meilleur $u$ est aussi le moins mauvais, on peut penser \`a choisir comme estimateur de $P\et$ la loi $P_{\widetilde{t_{n}}}$ avec
\begin{equation}
\widetilde{t}_{n}(\biX)=\argmin_{u\in\sbS}\left[\sup_{t\in\sbS}\sum_{i=1}^{n}\psi\left(\sqrt{\frac{t(X_{i})}{u(X_{i})}}\right)\right],
\label{eq-rhoest}
\end{equation}
en supposant ici, pour simplifier l'expos\'e, que l'$\argmin$ existe (mais n'est pas n\'ecessairement unique). Lorsque $\psi$ est la fonction $\log$ et $\widehat{t}_{n}$ l'e.m.v.,
\begin{align*}
\lefteqn{\argmin_{u\in\sbS}\left[\sup_{t\in\sbS}\sum_{i=1}^{n}\log\left(\sqrt{\frac{t(X_{i})}{u(X_{i})}}\right)\right]}\hspace{30mm}\\&=\frac{1}{2}
\argmin_{u\in\sbS}\left[\left(\sup_{t\in\sbS}\sum_{i=1}^{n}\log\left(t(X_{i})\st\right)\right)
-\sum_{i=1}^{n}\log\left(u(X_{i})\st\right)\right]\\&=\frac{1}{2}\argmin_{u\in\sbS}\left[\sum_{i=1}^{n}\log\left(\widehat{t}_{n}(X_{i})\st\right)-\sum_{i=1}^{n}\log\left(u(X_{i})\st\right)\right]\\&=\frac{1}{2}\argmax_{u\in\sbS}\left[\sum_{i=1}^{n}\log\left(u(X_{i})\st\right)-\sum_{i=1}^{n}\log\left(\widehat{t}_{n}(X_{i})\st\right)\right]\\&=
\argmax_{u\in\sbS}\sum_{i=1}^{n}\log\left(u(X_{i})\st\right)=\widehat{t}_{n}(\biX).
\end{align*}
On retrouve ainsi le maximum de vraisemblance parce que $\log(a/b)=\log a-\log b$, mais 
cet argument ne fonctionne plus avec une fonction $\psi\ne\alpha\log$ parce qu'alors
\[
\sum_{i=1}^{n}\psi\left(\sqrt{\frac{t(X_{i})}{u(X_{i})}}\right)\quad\mbox{ne s'\'ecrit pas}\quad
\sum_{i=1}^{n}\psi\left(\sqrt{t(X_{i})}\right)-\sum_{i=1}^{n}\psi\left(\sqrt{u(X_{i})}\right).
\]
Dans ce cas, on ne peut plus caract\'eriser $\widetilde{t}_{n}(\biX)$ d\'efini par (\ref{eq-rhoest}) comme 
\[
\argmax_{t\in\sbS}\left[\sum_{i=1}^{n}\psi\left(\sqrt{t(X_{i})}\right)\right].
\]

\section{Les fondamentaux des $\rho$-estimateurs\label{S4}}

\subsection{Cadre statistique et d\'efinitions\label{S4a}}
Comme nous voulons que ce qui va suivre s'applique aussi \`a des probl\`emes de r\'egression pour lesquels les observations ne sont pas n\'ecessairement i.i.d., nous allons travailler dans le cadre plus g\'en\'eral suivant, en omettant d\'elib\'er\'ement un certain nombre de points techniques (en particulier li\'es \`a des probl\`emes de mesurabilit\'e) qui ne pourraient que nuire \`a la clart\'e de cette pr\'esentation. 

On observe $\biX=(\etc{X})\in\sbX=\prod_{i=1}^{n}\sX_{i}$ et les v.a.\ $X_{i}\in\sX_{i}$ sont ind\'ependantes, de lois respectives $P\et_{i}$, donc le vecteur $\biX$ a pour loi $\gP\et=\bigotimes_{i=1}^{n}P\et_{i}$. Pour estimer $\gP\et$ \`a partir de $\biX$ on se donne un mod\`ele statistique $\sbP$ domin\'e par $\gmu=\bigotimes_{i=1}^{n}\mu_{i}$ avec $\gQ=\bigotimes_{i=1}^{n}Q_{i}$ pour tout $\gQ\in\sbP$, $Q_{i}=q_{i}\cdot\mu_{i}$ et $d\gQ/d\gmu=\gq=\bigotimes_{i=1}^{n}q_{i}$ $\gmu$-p.p., les versions des densit\'es pouvant \^etre choisies de mani\`ere arbitraire. On associe ainsi \`a $\sbP$ un ensemble $\sbQ$ de densit\'es $\gq$ par rapport \`a $\gmu$ et $\sbP=\{\gQ=\gq\cdot\gmu,\,\gq\in\sbQ\}$. Clairement, \'etant donn\'e $\sbP$, l'ensemble des fonctions $\gq$ n'est pas unique mais nous passerons ici sur cet aspect technique comme indiqu\'e pr\'ec\'edemment, les r\'esultats que nous pr\'esenterons ne d\'ependant pas du choix que l'on fait des versions de ces densit\'es. Il est important de noter que nous ne supposons pas que $\gP\et\in\sbP$ mais seulement qu'il s'agit d'une loi produit.

On se donne sur $[0,+\infty]$ une fonction $\psi$ v\'erifiant $\psi(1)=0$ et l'on consid\`ere 
la famille de statistiques
\begin{equation}
\gT(\gX,\gq,\gq')=\sum_{i=1}^{n}\psi\pa{\sqrt{{q_{i}'(X_{i})}\over {q_{i}(X_{i})}}}\,\mbox{ pour tous les }\:\gq,\gq'\in\sbQ.
\label{eq-defT}
\end{equation}
On d\'efinit alors un estimateur $\widehat{\gP}$ par
\begin{equation}
\widehat{\gP}=\widehat{\gp}\cdot\gmu=\bigotimes_{i=1}^{n}\left(\widehat{p}_{i}\cdot\mu_{i}\right)\quad\;\mbox{avec}\;\quad\widehat{\gp}=\argmin_{\gq\in\sbQ}\sup_{\gq'\in\sbQ}\gT(\gX,\gq,\gq').
\label{eq-rhoest}
\end{equation}
Pour simplifier, ici comme dans toute la suite, nous supposerons que l'$\argmin$ est atteint donc que $\widehat{\gp}$ existe. Notons qu'il n'est pas n\'ecessairement unique. Dans le cas g\'en\'eral, il suffit de choisir pour $\widehat{\gp}$ un \'el\'ement de $\sbQ$ qui minimise
$\sup_{\gq'\in\sbQ}\gT(\gX,\gq,\gq')$ approximativement.
%
\begin{rem}
Comme $\psi(1)=0$, $\gT(\gX,\gq,\gq)=0$, donc
\[
\min_{\gq\in\sbQ}\sup_{\gq'\in\sbQ}\gT(\gX,\gq,\gq')\ge0
\]
et si l'on trouve un \'el\'ement $\gq$ de $\sbQ$ tel que $\sup_{\gq'\in\sbQ}\gT(\gX,\gq,\gq')=0$, alors c'est un $\rho$-estimateur.
\end{rem}
%
\subsection{Choix de la fonction de perte\label{S4b}}

Pour \'evaluer les performances d'un estimateur de $\gP\et$, nous utiliserons des fonctions de perte li\'ees \`a la distance de Hellinger.
%
\begin{df}\label{def-Hell}
\'Etant donn\'e deux probabilit\'es $P$ et $Q$ et une mesure arbitraire $\nu$ qui domine $P$ et $Q$ on d\'efinit la distance de Hellinger $h$ et l'affinit\'e de Hellinger $\rho$ entre $P$ et $Q$ par
\[
h^{2}(P,Q)=\frac{1}{2}\int\left(\sqrt{\frac{dP}{d\nu}}-\sqrt{\frac{dP}{d\nu}}\right)^{2}d\nu
\quad\mbox{et}\quad\rho(P,Q)=\int\sqrt{\frac{dP}{d\nu}\frac{dQ}{d\nu}}d\nu.
\]
\end{df}
On notera que ces deux quantit\'es ne d\'ependent ni du choix de $\nu$ ni de celui des versions des deux densit\'es $dP/d\nu$ et $dQ/d\nu$ et que
\[
h^{2}(P,Q)=1-\rho(P,Q)\qquad\mbox{et}\qquad\rho\left(P^{\otimes n},Q^{\otimes n}\right)
=\rho^{n}(P,Q).
\]
Dans le cas des lois produits qui nous int\'eressent, soit $\gP=\bigotimes_{i=1}^{n}P_{i}$ et $\gQ=\bigotimes_{i=1}^{n}Q_{i}$, nous consid\'ererons la distance $H$ (\'egalement appel\'ee ici distance de Hellinger, bien que $H(\gP,\gQ)\ne h(\gP,\gQ)$) d\'efinie par
\[
H^{2}(\gP,\gQ)=\sum_{i=1}^{n}h^{2}(P_{i},Q_{i})\in[0,n].
\]
La qualit\'e d'un estimateur $\widehat{\gQ}(\biX)$ de la vraie loi $\gP\et$ sera mesur\'ee par son risque quadratique
\[
\E_{\gP\et}\left[H^{2}\left(\widehat{\gQ}(\biX),\gP\et\right)\right].
\]
Un cas particulier important est celui de lois produits correspondant \`a des v.a.\ i.i.d.~: $\gP=P^{\otimes n}$ et $\gQ=Q^{\otimes n}$ pour lequel $H^{2}(\gP,\gQ)=nh^{2}(P,Q)$.

\subsection{Quelles conditions sur la fonction $\psi$ ?\label{S4c}}
Comme nous l'avons expliqu\'e, l'id\'ee de cette construction est de remplacer, dans le calcul de la log-vraisemblance, la fonction $\log$ par une fonction $\psi$ born\'ee. N\'ean\-moins, afin 
de pr\'eserver certaines des propri\'et\'es du log, en particulier d'obtenir un analogue de (\ref{eq-test}), il convient que la fonction $\psi$ ait certaines similitudes avec la fonction $\log$. Nous demanderons donc que $\psi$ soit continue et strictement croissante avec $\psi(1/x)=-\psi(x)$. Comme elle est born\'ee et continue, la fonction $\psi$ v\'erifie
\[
\psi(+\infty)=\lim_{x\rightarrow+\infty}\psi(x)=\lim_{y\rightarrow0}\psi(1/y)=-\lim_{y\rightarrow0}\psi(y)=-\psi(0).
\]
Le remplacement de $\psi$ par $\alpha\psi$ avec $\alpha>0$ ne modifiant pas la valeur des
$\rho$-estimateurs, nous conviendrons de fixer $\psi(0)=-1$ donc $\psi(+\infty)=1$ et, comme nous voulons utiliser la distance de Hellinger pour mesurer les performances des $\rho$-estimateurs, nous imposerons les relations suivantes :

Quelles que soient les densit\'es $\gq,\gq'\in\sbQ$, la probabilit\'e $\gR\in\sbP$ et $i\in\{1,\ldots,n\}$,
\begin{equation}
\int_{\sX_{i}}\psi\pa{\sqrt{q'_{i}\over q_{i}}}\,dR_{i}\le a_{0}h^{2}(R_{i},Q_{i})-a_{1}h^{2}(R_{i},Q'_{i})
\label{eq-esp}
\end{equation}
et 
%
\begin{equation}
\int_{\sX_{i}}\psi^{2}\pa{\sqrt{q'_{i}\over q_{i}}}\,dR_i\le a_{2}^{2}\cro{h^{2}(R_{i},Q_{i})+h^{2}(R_{i},Q'_{i})},
\label{eq-var}
\end{equation}
o\`u les constantes $a_{0}$, $a_{1}$ et $a_{2}$ satisfont aux relations
\[
a_{0}\ge1\ge a_{1}\ge0\qquad\mbox{et}\qquad a_{2}^{2}\ge 1\vee(6a_{1}). 
\]
Pour que les rapports $q'_{i}/q_{i}$ soient toujours bien d\'efinis, nous utilisons les conventions $0/0=1$ et $a/0=+\infty$ pour tout $a>0$.

La relation (\ref{eq-esp}) montre que si $R_{i}$ est la vraie loi de $X_{i}$ et $Q_{i}$ est bien plus proche de $R_{i}$ que ne l'est $Q'_{i}$ alors l'int\'egrale de gauche qui est l'esp\'erance de $\psi\left(\sqrt{q'_{i}(X_{i})/q_{i}(X_{i})}\right)$ est n\'egative. Quant \`a (\ref{eq-var}) elle permet de contr\^oler la variance de la m\^eme quantit\'e. Comme la fonction $\psi$ est born\'ee, on peut en d\'eduire des in\'egalit\'es de d\'eviations de type Bernstein pour $\gT(\gX,\gq,\gq')$ qui est une somme de variables ind\'ependantes, born\'ees et dont on contr\^ole la variance. D'o\`u l'importance de cette bornitude. 

Dans la situation id\'eale o\`u l'on aurait $a_{0}=a_{1}=1$, on d\'eduirait de (\ref{eq-esp}) que
\[
\E_{\gP\et}\left[\gT(\gX,\gq,\gq')\right]=\sum_{i=1}^{n}\int_{\sX_{i}}\psi\pa{\sqrt{q'_{i}\over q_{i}}}\,dP_i\et\le H^{2}(\gP\et,\gQ)-H^{2}(\gP\et,\gQ')
\]
et, en inversant les r\^oles de $\gQ$ et $\gQ'$, que
\begin{align*}
-\sum_{i=1}^{n}\int_{\sX_{i}}\psi\pa{\sqrt{q'_{i}\over q_{i}}}\,dP_i\et&=\sum_{i=1}^{n}\int_{\sX_{i}}\psi\pa{\sqrt{q_{i}\over q'_{i}}}\,dP_i\et
=\E_{\gP\et}\left[\gT(\gX,\gq',\gq)\right]\\&\le H^{2}(\gP\et,\gQ')-H^{2}(\gP\et,\gQ),
\end{align*}
parce que $\psi(1/x)=-\psi(x)$. Finalement, comme $\left[\gT(\gX,\gq',\gq)\right]=-\left[\gT(\gX,\gq,\gq')\right]$,
\[
-\left[H^{2}(\gP\et,\gQ')-H^{2}(\gP\et,\gQ)\right]\le\E_{\gP\et}\left[\gT(\gX,\gq,\gq')\right]\le H^{2}(\gP\et,\gQ)-H^{2}(\gP\et,\gQ'),
\]
donc
\[
\E_{\gP\et}\left[\gT(\gX,\gq,\gq')\right]=H^{2}(\gP\et,\gQ)-H^{2}(\gP\et,\gQ').
\]
De m\^eme 
\[
\Var_{\gP\et}\left(\st\gT(\gX,\gq,\gq')\right)\le a_{2}\left[H^{2}(\gP\et,\gQ)+H^{2}(\gP\et,\gQ')\right].
\]
La premi\`ere \'egalit\'e est \`a rapprocher, dans le cadre i.i.d., de
\begin{align*}
\lefteqn{\E_{\gP\et}\left[\sum_{i=1}^{n}\log\left(\sqrt{\frac{q'(X_{i})}{q(X_{i})}}\right)\right]}\hspace{20mm}\\&=
\E_{\gP\et}\left[\sum_{i=1}^{n}\log\left(\sqrt{\frac{p\et(X_{i})}{q(X_{i})}}\right)\right]-\E_{\gP\et}\left[\sum_{i=1}^{n}\log\left(\sqrt{\frac{p\et(X_{i})}{q'(X_{i})}}\right)\right]\\&=(n/2)\cro{K(P\et,Q)-K(P\et,Q')},
\end{align*}
pourvu que les deux divergences de Kullback-Leibler, $K(P\et,Q)$ et $K(P\et,Q')$, soient finies, ce qui n'est absolument pas garanti, $P\et$ \'etant inconnue. 

Nous ignorons combien de fonctions satisfont aux relations (\ref{eq-esp}) et (\ref{eq-var}) mais nous n'en connaissons que deux :
\[
\psi_{1}(x)={x-1\over x+1}\qquad\mbox{et}\qquad\psi_{2}(x)={x-1\over \sqrt{x^{2}+1}},
\]
la premi\`ere fonction $\psi_{1}(x)$ \'etant la plus simple et tr\`es proche de la fonction $x\mapsto(1/2)\log(x)$ au voisinage de 1 puisque
\[
0,96<\frac{\psi_{1}(x)}{(\log x)/2}\le1\quad\mbox{pour }\;1/2\le x\le2
\]
et
\[
0,86<\frac{\psi_{1}(x)}{(\log x)/2}\le1\quad\mbox{pour }\;1/4\le x\le4.
\]
Pour la fonction $\psi_{1}$, $a_{0}=4$, $a_{1}=3/8$ et $a_{2}^{2}=3\sqrt{2}$, et pour la fonction $\psi_{2}$, $a_{0}=4,97$, $a_{1}=0,083$ et $a_{2}^{2}=3+2\sqrt{2}$.

\section{Les performances des $\rho$-estimateurs\label{S5}}

\subsection{Un r\'esultat g\'en\'eral\label{S5a}}
La qualit\'e d'un $\rho$-estimateur $\widehat{\gP}$ d\'efini par (\ref{eq-rhoest}), mesur\'ee par son risque quadratique $\E_{\gP\et}\left[H^{2}\left(\widehat{\gP}(\biX),\gP\et\right)\right]$, ne
d\'epend que du mod\`ele $\sbP$ choisi et de la vraie loi $\gP\et$. Plus pr\'ecis\'ement,
%
\begin{thm}\label{thm-Main}
Pour toute mesure-produit $\gP\et=\bigotimes_{i=1}^{n}P\et_{i}$, toute probabilit\'e $\overline \gP\in\sbP$ et tout r\'eel positif $\xi$,
\begin{equation}
\P_{\gP\et}\cro{H^{2}(\gP\et,\widehat \gP)\le C\left[H^{2}(\gP\et,\overline\gP)+D_{n}^{\sbP}(\gP\et,\overline \gP)+\xi\right]}\ge 1-e^{-\xi},
\label{eq-1}
\end{equation}
o\`u $C\ge 1$ est une constante num\'erique universelle (on peut prendre $C=26350$) et $D_{n}^{\sbP}(\gP\et,\overline \gP)\in [1,n]$ un terme de dimension locale (relative \`a $\overline{\gP}$) du mod\`ele $\sbP$.
\end{thm}
En int\'egrant par rapport \`a $\xi$ l'in\'egalit\'e pr\'ec\'edente, on obtient, puisque
$D_{n}^{\sbP}(\gP\et,\overline \gP)\ge1$,
\begin{align}
\E_{\gP\et}\cro{H^{2}(\gP\et,\widehat \gP)}&\le C\inf_{\overline \gP\in\sbP}\left\{H^{2}(\gP\et,\overline\gP)+D_{n}^{\sbP}(\gP\et,\overline \gP)+1\right\}\nonumber\\&\le
C\inf_{\overline \gP\in\sbP}\left\{H^{2}(\gP\et,\overline\gP)+2D_{n}^{\sbP}(\gP\et,\overline \gP)\right\}.
\label{eq-1b}
\end{align}
Dans la suite, $C$ d\'esignera une constante num\'erique qui pourra changer de ligne en ligne. La dimension $D_{n}^{\sbP}(\gP\et,\overline \gP)$ est d\'efinie \`a partir des fluctuations du processus $\gT(\gX,\overline{\gp},\gq)$ pour $\gq\in\sbQ$ autour de son esp\'erance sous $\gP\et$. Dans un grand nombre de situations (mais pas toujours), on peut borner cette dimension uniform\'ement en $\overline{\gP}$ et $\gP\et$ par une quantit\'e $D_{n}(\sbP)$ qui ne d\'epend que du mod\`ele que l'on utilise. Dans ce cas, la borne (\ref{eq-1b}) devient
\begin{equation}
\E_{\gP\et}\cro{H^{2}(\gP\et,\widehat \gP)}\le C\left[H^{2}(\gP\et,\sbP)+D_{n}(\sbP)\right],
\label{eq-1c}
\end{equation}
avec
\[
H^{2}(\gP\et,\sbP)=\inf_{\overline \gP\in\sbP}H^{2}(\gP\et,\overline\gP).
\]
En particulier,
\begin{equation}
\sup_{\gP\et\in\sbP}E_{\gP\et}\cro{H^{2}(\gP\et,\widehat \gP)}\le CD_{n}(\sbP).
\label{eq-rismax}
\end{equation}
On retrouve dans (\ref{eq-1c}) une majoration classique du risque par la somme d'un terme d'approximation $H^{2}(\gP\et,\sbP)$ et d'un terme de dimension $D_{n}(\sbP)$, lequel caract\'erise la ``taille" du mod\`ele $\sbP$. Cette formule met clairement en \'evidence la robustesse d'un tel estimateur \`a un \'ecart possible de la vraie loi au mod\`ele : si $\gP\et$ appartient au mod\`ele $\sbP$ le risque est major\'e par la dimension de $\sbP$ qui permet de contr\^oler le risque maximum sur $\sbP$ par (\ref{eq-rismax}), sinon il faut y ajouter un terme d'approximation proportionnel \`a $H^{2}(\gP\et,\sbP)$ traduisant le fait que le mod\`ele $\sbP$ est inexact. Si l'erreur de mod\'elisation est suffisamment faible, ce terme sera petit et la borne de risque comparable \`a celle que l'on obtient lorsque le mod\`ele est exact.

\subsection{Quelques cas particuliers\label{S5b}}
--- Dans le cadre de l'estimation d'une densit\'e \`a partir de $n$ observations i.i.d.\ de densit\'e $p\et$ avec $\gP\et=(p\et\cdot\mu)^{\otimes n}$ et d'un mod\`ele $\sbP=\{P^{\otimes n}, P\in\sP\}$, les $\rho$-estimateurs prennent la forme $\widehat \gP=\widehat{P}^{\otimes n}$ avec $\widehat{P}\in\sP$ et la borne de risque (\ref{eq-1c}) s'\'ecrit alors
\begin{equation}\label{eq-RisqueMax}
\E_{\gP\et}\cro{h^{2}(P\et,\widehat P)}\le C\left[h^{2}(P\et,\sP)+n^{-1}D_{n}(\sbP)\right].
\end{equation}

--- Lorsque le mod\`ele $(\sbP,H)$ est de dimension m\'etrique finie major\'ee par $\overline{D}$, ce qui signifie que l'on peut recouvrir toute boule de $\sbP$ de rayon $x\eta$ (avec $x\ge2$) par un nombre de boules de rayon $\eta$ major\'e par $\exp\left(\overline{D}x^{2}\right)$, alors 
$D_{n}^{\sbP}(\gP\et,\overline \gP)\le C\overline{D}$. Dans ce cas, les $\rho$-estimateurs ont les m\^emes performances que les T-estimateurs de Birg{\'e}~\citeyearpar{MR2219712} mais, comme dans le cas des T-estimateurs, l'hypoth\`ese de dimension m\'etrique finie  n\'ecessite que l'espace m\'etrique $(\sbP,H)$ soit compact. 

--- Dans un cadre de  r\'egression \`a ``plan fix\'e" pour lequel $X_{i}=(w_{i},Y_{i})$, $\lambda$ est la mesure de Lebesgue sur $\R$ et 
\[
Y_{i}=f\et(w_{i})+\varepsilon_{i}\quad\mbox{avec }\varepsilon_{i}\mbox{ de loi }s\cdot\lambda,
\]
on utilise un mod\`ele correspondant aux lois possibles des $Y_{i}$ donn\'e par
\[
Y_{i}=f(w_{i})+\varepsilon_{i}\quad\mbox{avec }\varepsilon_{i}\mbox{ de loi }r\cdot\lambda,
\]
o\`u $r$ est une densit\'e fix\'ee et $f\in \cF$. Si la famille $\cF$ est incluse dans un espace vectoriel de dimension $k$ et la densit\'e $r$ est unimodale, on peut montrer que $D_{n}^{\sbP}(\gP\et,\overline \gP)\le Ck\log n$ en utilisant des arguments de classes de Vapnik. On perd alors un facteur $\log n$ dans les bornes de risque mais on s'est d\'ebarass\'e des hypoth\`eses de compacit\'e puisque $\cF$ peut \^etre un espace vectoriel. Nous ignorons si un tel facteur logarithmique est n\'ecessaire ou non si l'on se restreint \`a des estimateurs robustes au sens ci-dessus.

\section{Liens entre $\rho$-estimateurs et e.m.v.}
Consid\'erons un mod\'ele param\'etrique tr\`es simple dans lequel les $X_{i}$ sont suppos\'ees i.i.d.\  et uniformes sur $[\theta-1,\theta+1]$. On peut montrer dans ce cas que, si le mod\`ele est exact, c'est-\`a-dire si les $X_{i}$ sont effectivement i.i.d.\  et uniformes sur $[\theta\et-1,\theta\et+1]$ pour un $\theta\et\in\R$, alors l'estimateur $\left(X_{(1)}+X_{(n)}\right)\!/2$, qui est un estimateur du maximum de vraisemblance, est aussi un $\rho$-estimateur. La diff\'erence entre l'e.m.v.\ et le $\rho$-estimateur est que, si le mod\`ele est l\'eg\`erement inexact, le $\rho$-estimateur donnera un r\'esultat sens\'e, c'est-\`a-dire une estimation par une loi uniforme sur $[\widetilde{\theta}-1,\widetilde{\theta}+1]$ proche de la vraie loi, alors que le maximum de vraisemblance peut fort bien ne pas exister, par exemple si $X_{(n)}-X_{(1)}>2$.

Dans un mod\`ele param\'etrique tr\`es r\'egulier $\{p_{\theta}, \theta\in[0,1]\}$ avec des observations i.i.d., si les estimateurs du maximum de vraisemblance $p_{\widehat{\theta}_{n}}$ sont consistants alors ce sont aussi des $\rho$-estimateurs avec une probabilit\'e tendant vers 1 lorsque $n\rightarrow+\infty$ --- cf.\  le th\'eor\`eme~19 de Baraud, Birg\'e and Sart~\citeyearpar{MR3595933} ---.

Toujours dans un cadre i.i.d.,\ si le mod\`ele consiste en un ensemble convexe $\sQ$ de densit\'es par rapport \`a une mesure de r\'ef\'erence $\mu$, et que la fonction $\psi$ vaut $\psi_{1}$ ou $\psi_{2}$, alors l'application
\[
(u,t)\mapsto \sum_{i=1}^{n}\psi\left(\sqrt{\frac{t(X_{i})}{u(X_{i})}}\right)
\]
admet un point selle sur $\sQ\times \sQ$ car l'application $(x,y)\mapsto \psi(\sqrt{x/y})=-\psi(\sqrt{y/x})$ est concave en $x$ pour $y>0$ fix\'e et convexe en $y$ pour $x>0$ fix\'e. Ce point selle est \`a la fois un $\rho$-estimateur et l'e.m.v.\ sur $\sQ$. Ce r\'esultat, qui repose sur un argument d\^u \`a Su Weijie (2016, communication personnelle), est d\'emontr\'e dans la section~6 de Baraud et Birg\'e~\citeyearpar{Baraud:2016kq}. Une telle situation, pour laquelle l'e.m.v.\ est un $\rho$-estimateur et en a donc toutes les propri\'et\'es, se retrouve dans les exemples suivants: 

\begin{lista}
\item  l'ensemble des densit\'es sur $\sX$ qui sont constantes par morceaux sur une partition fix\'ee de $\sX$ de cardinal $D$. Dans ce cas, notons que l'e.m.v.\ est alors l'histogramme associ\'e \`a cette partition;
\item  l'ensemble des densit\'es d\'ecroissantes sur $[0,+\infty)$ ou croissantes sur $(-\infty,0]$ ou unimodales sur $\R$ avec un mode en 0. En particulier, l'estimateur de Grenander pour les densit\'es d\'ecroissantes sur $\R_{+}$, qui est l'e.m.v., est un $\rho$-estimateur;
\item l'enveloppe convexe d'une nombre fini de densit\'es (estimateurs pr\'eliminaires obtenus \`a partir d'un \'echantillon ind\'ependant), ce qui nous fournit le cadre de l'agr\'egation convexe.
\end{lista}
Un r\'esultat analogue pour des cadres statistiques plus g\'en\'eraux que l'estimation de densit\'e a \'et\'e montr\'e par Mathieu Sart (2017, communication personnelle).

En r\'esum\'e, sous des conditions convenables, le maximum de vraisemblance est un $\rho$-estimateur, exactement (cadre convexe) ou asymptotiquement (mod\`ele param\'etrique suffisamment r\'egulier).

En revanche, si la vraie loi d\'evie un peu du mod\`ele, il se peut, lorsque l'on ne se trouve pas dans le cadre de variables i.i.d.\ et d'un mod\`ele convexe, que le maximum de vraisemblance en soit gravement affect\'e alors que les $\rho$-estimateurs ne le seront que faiblement.

Pour en revenir \`a l'exemple de la section~\ref{sect-EMVbad}, nous allons v\'erifier que le $\rho$-estimateur sur le mod\`ele associ\'e \`a l'espace des param\`etres $\Theta'$ est l'estimateur du maximum de vraisemblance $\widetilde \theta=(X_{0},0)$. Donc que son risque quadratique est born\'e par 5.

En choisissant comme mesure de r\'eference la mesure $\gmu$ ayant pour densit\'e (par rapport \`a la mesure de Lebesgue sur $\R^{k+1}$) la fonction
\[
(x_{0},x')\mapsto(2\pi)^{-(k+1)/2}\exp\left[-\left(x_{0}^{2}/2\right)+\left(\|x'\|^{2}/2\right)\right],
\]
les lois du mod\`ele statistique $\{\gP_{\theta},\; \theta\in \Theta'\}=\{\gP_{(\theta_{0},0)},\; \theta_{0}\in \R\}$ ont pour densit\'es  respectives (par rapport \`a $\gmu$) $\gp_{\theta_{0}}:(x_{0},x')\mapsto \exp\cro{-x_{0}\theta_{0}+(\theta_{0}^{2}/2)}$ de sorte que le $\rho$-estimateur $\widehat \theta_{0}$ du param\`etre $\theta_{0}$, au vu de la seule observation $\bsX=(X_{0},\ldots,X_{k})$, ce qui correspond \`a $n=1$ dans (\ref{eq-defT}), minimise l'application 
\[
\theta_{0}\mapsto\sup_{\theta_{0}'\in\R}\psi\pa{\sqrt{\gp_{\theta_{0}'}(\bsX)\over \gp_{\theta_{0}}(\bsX)}}=\sup_{\theta_{0}'\in\R}\psi\pa{\frac{1}{4}\exp\left[-2X_{0}\theta_{0}'+(\theta_{0}')^{2}-\st\pa{-2X_{0}\theta_{0}+\theta_{0}^{2}}\right]}.
\]
La fonction $\psi$ \'etant strictement croissante, le maximum, pour un $\theta_{0}$ fix\'e, est atteint lorsque $\theta'_{0}=X_{0}$ (ind\'ependamment de $\theta_{0}$) et le minimum en $\theta_{0}$ atteint au point $\widehat \theta_{0}=X_{0}$. Le $\rho$-estimateur du param\`etre $\theta=(\theta_{0},\theta')$ est donc l'estimateur du maximum de vraisemblance $\widetilde \theta$.

\section{Extension\label{S7}}
Afin de rendre cette pr\'esentation aussi simple que possible, nous nous sommes content\'es jusqu'ici d'exposer la construction des $\rho$-estimateurs sur un seul mod\`ele $\sbP$ mais, comme dans le cadre des T-estimateurs ou autres, il est possible de travailler avec une famille d\'enombrable $\{\sbP_{m}, m\in\cM\}$ de mod\`eles simultan\'ement en ajoutant \`a la statistique $\gT$ une p\'enalit\'e. 

Comme pr\'ec\'edemment, chaque mod\`ele $\sbP_{m}$ est d\'ecrit par un ensemble $\sbQ_{m}$ de densit\'es par rapport \`a $\gmu$, c'est-\`a-dire que $\sbP_{m}=\{\gQ=\gq\cdot\gmu,\,\gq\in\sbQ_{m}\}$ pour tout $m$. En outre, $\sbP_{m}$ est affect\'e d'un poids $\exp[-\Delta(m)]$ avec
\begin{equation}
\sum_{m\in\cM}e^{-\Delta(m)}=1,
\label{eq-delta}
\end{equation}
de sorte que cette famille de poids peut \^etre consid\'er\'ee comme une loi a priori sur l'ensemble des mod\`eles. Pour simplifier la pr\'esentation et bien que ceci ne soit en aucun cas n\'ecessaire, nous supposerons que tous nos mod\`eles sont disjoints de sorte qu'\`a chaque loi $\gQ\in\sbP=\bigcup_{m\in\cM}\sbP_{m}$ correspond un unique $m$ tel que $\gQ\in\sbP_{m}$. On associe alors \`a chaque $\gQ\in\sbP$ une p\'enalit\'e $\pen(\gQ)\ge0$ telle que, pour tout $\overline{\gP}\in\sbP$,
\begin{equation}
G(\gP\et,\overline\gP)+\pen(\gQ)\ge\kappa\cro{D^{\sbP_{m}}_{n}(\gP\et,\overline \gP)+\Delta(m)}
\,\mbox{ pour tout }\gQ\in\sbP_{m},
\label{eq-L}
\end{equation}
o\`u $\kappa$ est une constante num\'erique et $G(\gP\et,\cdot)$ une fonction sur $\sbP$. Notons que si l'on ne dispose que d'un seul mod\`ele $\sbP_{0}$, c'est-\`a-dire si $\cM=\{0\}$, on peut fixer $\Delta(0)=0=\pen(\gQ)$ pour tout $\gQ\in\sbP_{0}$ et $G(\gP\et,\overline\gP)=\kappa D^{\sbP_{0}}_{n}(\gP\et,\overline \gP)$.

Dans la situation type o\`u $D^{\sbP_{m}}_{n}(\gP\et,\overline \gP)=D_{n}(\sbP_{m})$ ne d\'epend que du mod\`ele $\sbP_{m}$, on peut prendre la fonction $G$ identiquement nulle et 
\begin{equation}
\pen(\gQ)=\kappa\cro{D_{n}(\sbP_{m})+\Delta(m)}\,\mbox{ pour tout }\gQ\in\sbP_{m},
\label{eq-penQ}
\end{equation}
c'est-\`a-dire une p\'enalit\'e constante sur chaque mod\`ele. 

Dans ce cadre, la d\'efinition des $\rho$-estimateurs doit \^etre modifi\'ee comme suit. L'on d\'efinit
\begin{equation}
\gup(\bsX,\gq)=\sup_{\gq'\in\sbQ}\cro{\gT(\bsX,\gq,\gq')-\pen(\gQ')}+\pen(\gQ)\,\mbox{ pour tout }\gq\in\sbQ
\label{eq-gup}
\end{equation}
et un $\rho$-estimateur s'\'ecrit $\widehat{\gP}=\widehat{\gp}\cdot\gmu$ avec $\widehat{\gp}=\argmin_{\gq\in\sbQ}\gup(\bsX,\gq)$. Tout $\rho$-estimateur $\widehat{\gP}$ satisfait alors au r\'esultat suivant.
%
\begin{thm}\label{thm-Mainsel}
Pour toute mesure-produit $\gP\et=\bigotimes_{i=1}^{n}P\et_{i}$, toute probabilit\'e $\overline \gP\in\sbP$ et tout r\'eel positif $\xi$,
\begin{equation}
\P_{\gP\et}\cro{H^{2}(\gP\et,\widehat \gP)\le \overline{C}\left[H^{2}(\gP\et,\overline\gP)+G(\gP\et,\overline \gP)+\pen(\overline{\gP})+\xi\right]}\ge 1-e^{-\xi},
\label{eq-2}
\end{equation}
o\`u $\overline{C}=\overline{C}(a_{0},a_{1},a_{2})\ge 1$ d\'epend uniquemement de $a_{0}$, $a_{1}$ et $a_{2}$.
\end{thm}
Dans la suite, comme ici, nous d\'esignerons par $\overline{C}$ des quantit\'es qui peuvent d\'ependre de certains param\`etres correspondant aux hypoth\`eses que nous ferons, afin de les distinguer des constantes universelles not\'ees $C$ ou $C'$. Lorsque la fonction $G$ est nulle et la p\'enalit\'e est donn\'ee par (\ref{eq-penQ}), (\ref{eq-2}) devient
\[
\P_{\gP\et}\!\cro{H^{2}(\gP\et,\widehat \gP)\le \overline{C}\inf_{m\in\cM}\left[H^{2}(\gP\et,\sbP_{m})+D_{n}(\sbP_{m})+\Delta(m)+\xi\right]}\ge 1-e^{-\xi},
\]
et, apr\`es int\'egration par rapport \`a $\xi>0$,
\begin{equation}
\E_{\gP\et}\cro{H^{2}(\gP\et,\widehat \gP)}\le \overline{C}\inf_{m\in\cM}\left[H^{2}(\gP\et,\sbP_{m})+D_{n}(\sbP_{m})+\Delta(m)\right].
\label{eq-1cb}
\end{equation}
Notons que si l'on utilisait le seul mod\`ele $m$, on d\'eduirait de (\ref{eq-1c}) que
\begin{equation}
\E_{\gP\et}\cro{H^{2}(\gP\et,\widehat \gP)}\le C\left[H^{2}(\gP\et,\sbP_{m})+D_{n}(\sbP_{m})\right],
\label{eq-1f}
\end{equation}
et l'in\'egalit\'e~\eref{eq-1cb} conduirait alors \`a un analogue de (\ref{eq-1c}) avec un terme suppl\'ementaire $\Delta(m)$. Si $\Delta(m)$ est au plus du m\^eme ordre de grandeur que $D_{n}(\sbP_{m})$ pour tout $m\in\cM$, on obtient l'exact analogue de (\ref{eq-1c}) \`a la constante $\overline{C}$ pr\`es et la s\'election de mod\`ele ne co\^ute rien au sens o\`u la borne de risque obtenue n'est jamais beaucoup plus grande que celle que l'on obtiendrait en faisant un choix a priori d'un mod\`ele dans la famille. Cette borne est en fait du m\^eme ordre de grandeur que celle que l'on obtiendrait en utilisant seulement le ``meilleur" mod\`ele, c'est-\`a-dire celui qui minimise en $m$ la quantit\'e
$H^{2}(\gP\et,\sbP_{m})+D_{n}(\sbP_{m})$ et optimise ainsi la borne (\ref{eq-1f}).

\section{Introduction aux propri\'et\'es des $\rho$-esti\-mateurs\label{S8}}
Consid\'erons d'abord ici la situation d'un seul mod\`ele $\sbP$ pour lequel la borne (\ref{eq-1b}) est valide, c'est-\`a-dire que
\[
\E_{\gP\et}\cro{H^{2}(\gP\et,\widehat \gP)}\le C\inf_{\overline \gP\in\sbP}\left\{H^{2}(\gP\et,\overline\gP)+D_{n}^{\sbP}(\gP\et,\overline \gP)\right\}.
\]
La quantit\'e $D_{n}^{\sbP}(\gP\et,\overline \gP)$ d\'epend en th\'eorie de la loi $\gP\et$ des observations mais dans toutes les situations que nous avons \'etudi\'ees il est possible de la majorer par une quantit\'e $D_{n}^{\sbP}(\overline \gP)$ qui ne d\'epend que du mod\`ele $\sbP$ et de $\overline \gP$. Dans ce cas, notre borne de risque devient
\begin{equation}
\E_{\gP\et}\cro{H^{2}(\gP\et,\widehat \gP)}\le C\left[H^{2}(\gP\et,\overline\gP)+D_{n}^{\sbP}(\overline \gP)\right]
\,\mbox{ pour tout }\overline{\gP}\in\sbP.
\label{eq-3b}
\end{equation}
Lorsque $\gP\et$ est un \'el\'ement $\overline \gP$ de $\sbP$ la borne devient 
\[
\E_{\gP\et}\cro{H^{2}({\gP\et},\widehat \gP)}=\E_{\overline{\gP}}\cro{H^{2}(\overline{\gP},\widehat \gP)}\le CD_{n}^{\sbP}(\overline{\gP}).
\]
D'apr\`es \eref{eq-3b}, si  $\gP\et\ne\overline{\gP}$, la borne ne se d\'et\'eriore pas plus que d'une quantit\'e $CH^{2}(\gP\et,\overline\gP)$, ce qui signifie en particulier que la borne de risque reste stable (comme fonction de $\gP\et$) au voisinage de chaque point $\overline\gP$ du mod\`ele. Il s'agit l\`a de la propri\'et\'e fondamentale des $\rho$-estimateurs que nous allons \`a pr\'esent illustrer. Pour ce faire, il conviendra de bien distinguer la vraie loi $\gP\et$ des observations $\biX$ du mod\`ele statistique $\sbP$ que nous introduisons pour construire nos estimateurs. La seule hypoth\`ese que nous faisons sur $\gP\et$ est que c'est une loi produit (qui n'appartient pas n\'ecessairement au mod\`ele). Le choix du mod\`ele est fond\'e sur certaines hypoth\`eses dont on esp\`ere seulement qu'elles ne sont pas trop erron\'ees si l'on veut que $\gP\et$ ne soit pas trop loin de $\sbP$.

Tous les mod\`eles statistiques que nous consid\'ererons dans la suite supposent les observations i.i.d.\ (alors qu'elles ne le sont pas forc\'ement). Cela revient \`a choisir un mod\`ele statistique $\sP$ pour la loi marginale des observations de sorte que 
\[
\overline \gP=\overline P^{\otimes n}\,\mbox{ avec }\,\overline P \in \sP\quad \mbox{pour tout }\overline \gP\in\sbP.
\]
Un $\rho$-estimateur s'\'ecrit donc $\widehat \gP=\widehat P^{\otimes n}$ avec $\widehat P\in\sP$ alors que $\gP\et=\bigotimes_{i=1}^{n}P_{i}\et$.

\section{Robustesse\label{S9}}

\subsection{Robustesse \`a l'hypoth\`ese d'\'equidistribution\label{S9a}}
Nous supposons ici qu'il existe une loi marginale $\overline P\in\sP$ telle que $P_{i}\et=\overline P$ pour tous les indices $i\in \{1,\ldots,n\}\setminus I$ o\`u $I$ est un sous-ensemble de $\{1,\ldots,n\}$ et si $i\in I$, la loi $P_{i}\et$ est arbitraire. Par exemple, si pour tout $i\in I$, $P_{i}\et=\delta_{x_{i}}$ (mesure de Dirac au point $x_{i}$) o\`u les $x_{i}$ sont des points arbitraires, cela signifie que notre suppos\'e $n$-\'echantillon contient en fait $|I|$ valeurs atypiques (o\`u $|I|$ d\'esigne le cardinal de l'ensemble $I$). Nous pouvons alors \'ecrire que 
\[
H^{2}(\gP\et,\widehat \gP)\ge \sum_{i\not \in I}h^{2}(P_{i}\et,\widehat P)=(n-|I|)h^{2}(\overline P,\widehat P)
\]
et 
\[
H^{2}(\gP\et,\overline  \gP)= \sum_{i\not \in I}h^{2}(P_{i}\et,\overline P)+\sum_{i\in I}h^{2}(P_{i}\et,\overline P)\le |I|.
\]
Il d\'ecoule donc de~\eref{eq-3b} que 
\[
\E_{\gP\et}\cro{h^{2}(\overline P,\widehat P)}\le C{|I|+D_{n}^{\sbP}(\overline \gP)\over n-|I|}.
\]

\begin{lista}
\item Lorsque $I=\varnothing$, c'est-\`a-dire lorsque les observations sont vraiment i.i.d. et que le mod\`ele est exact, on obtient la borne $CD_{n}^{\sbP}(\overline \gP)/n$ comme attendu.
\item  Si les donn\'ees ne sont pas exactement i.i.d., c'est-\`a-dire si $I\neq \varnothing$, mais que son cardinal n'est pas trop grand devant $D_{n}^{\sbP}(\overline \gP)\le n$, nous obtenons une borne de risque du m\^eme ordre de grandeur que la pr\'ec\'edente. En parti\-culier, le risque de l'estimateur reste stable \`a un possible \'ecart \`a l'hypoth\`ese d'\'equidistribution et notamment \`a la pr\'esence de quelques valeurs atypiques. 
\end{lista}

\subsection{Robustesse \`a la contamination\label{S9b}}
Supposons \`a pr\'esent que les donn\'ees sont vraiment i.i.d.\ et qu'il existe une loi $\overline P\in\sP$ et un nombre $\eps\in [0,1]$ (typiquement petit) tels que  
\[
P_{i}\et=P\et=(1-\eps)\overline P+\eps Q\quad \mbox{pour tout $i\in \{1,\ldots,n\}$},
\]
o\`u $Q$ est une loi arbitraire. Cela revient \`a supposer qu'une proportion $1-\eps$ de l'\'echantillon est correctement mod\'elis\'ee par une loi $\overline P$ du mod\`ele, qu'une proportion $\eps$ est distribu\'ee selon une loi $Q$ quelconque et que ce dernier \'echantillon vient contaminer le premier.  

Comme $h^{2}(P,R)\le\|P-R\|$ o\`u $\|P-R\|$ d\'esigne la distance en variation totale,
\[
h^{2}\left(P\et,\overline P\right)=h^{2}\left((1-\eps)\overline P+\eps Q,\overline P\right)\le \norm{\eps (Q-\overline P)}\le \eps,
\]
de sorte que 
\[
H^{2}(\gP\et,\overline\gP)=nh^{2}(P\et,\overline P)\le n\varepsilon
\]
et~\eref{eq-3b} conduit alors \`a la borne de risque 
\[
\E_{\gP\et}\cro{h^{2}(P\et ,\widehat P)}\le C\left[\eps+{D_{n}^{\sbP}(\overline \gP)\over n}\right].
\]
Cette borne reste donc stable \`a une contamination possible des donn\'ees dans une proportion $\eps$, tant que celle-ci n'est pas trop grande par rapport \`a $D_{n}^{\sbP}(\overline \gP)/n$. 

\section{Deux exemples\label{S10}}
Nous allons \`a pr\'esent consid\'erer deux types de mod\`eles statistiques particuliers et \'etudier les propri\'et\'es des $\rho$-estimateurs pour chacun d'eux. Nous ne reviendrons plus dans la suite sur les propri\'et\'es de robustesse vues ci-dessus.

\subsection{Mod\`eles de densit\'es sous contrainte de forme\label{S10a}}
Nous supposerons ici que les variables $X_{i}$ sont i.i.d.\  de densit\'e $s\et$ par rapport \`a la mesure de Lebesgue $\lambda$ sur $\R$ de sorte que $\gP\et=\gP_{\!s\et}=(s\et\cdot \lambda)^{\otimes n}$ et
\[
H^{2}(\gP_{\!s\et},\gP_{\!s})=nh^{2}(s\et\cdot \lambda,s\cdot \lambda)\quad \mbox{avec $\gP_{\!s}=(s\cdot \lambda)^{\otimes n}$.}
\] 
Pour simplifier, nous noterons $h(s\et,s)$ pour $h(s\et\cdot \lambda,s\cdot \lambda)$.

Le mod\`ele statistique consiste \`a supposer que $s\et$ appartient \`a une famille $S$ de densit\'es  d\'efinies par des propri\'et\'es qualitatives. Par exemple, $S$ est l'ensemble des densit\'es d\'ecroissantes sur $[0,+\infty)$. Nous allons montrer comment on peut analyser les performances du $\rho$-estimateur de $s\et$ sur ce type de mod\`eles qui sont en g\'en\'eral tr\`es gros (non-compacts pour la distance $H$) et pour lesquels il n'existe pas de vitesse d'estimation minimax. 

Nous n'allons en fait pas traiter l'exemple typique des densit\'es d\'ecroissantes sur $[0,+\infty)$ puisque, comme nous l'avons dit pr\'ec\'edemment, Su Weijie a d\'emontr\'e que dans cette situation le $\rho$-estimateur co\"incide exactement avec l'estimateur du maximum de vraisemblance, c'est-\`a-dire l'estimateur de Grenander, bien connu et tr\`es abondamment \'etudi\'e --- cf.\ par exemple Birg{\'e}~\citeyearpar{MR1026298} et les r\'ef\'erences incluses dans cet article ---.

Nous allons plut\^ot nous int\'eresser  \`a un probl\`eme plus complexe en choisissant pour $S$ l'ensemble de toutes les densit\'es qui sont monotones sur une demi-droite et nulles ailleurs. Cet ensemble contient les densit\'es pr\'ec\'edentes et plus g\'en\'eralement toutes les densit\'es d\'ecroissantes sur un intervalle de la forme $[a,+\infty)$, mais aussi toutes celles qui sont croissantes sur un intervalle  de la forme $(-\infty,a]$ ainsi que toutes les densit\'es des lois uniformes sur un intervalle compact. L'ensemble $S$ est stable par translation et changement d'\'echelle. Notre mod\`ele statistique est donc 
\[
\sbP=\{\gP_{\!s}=(s\cdot \lambda)^{\otimes n},\ s\in S\}.
\]
Bien que ce mod\`ele soit tr\`es riche (il n'existe pas de vitesse minimax sur un tel ensemble), le $\rho$-estimateur $\widehat \gP=\gP_{\widehat s}$ sur $\sbP$ n'est pas d\'eg\'en\'er\'e, comme nous allons le voir, et nous allons pouvoir \'etudier ses propri\'et\'es \`a partir de l'in\'egalit\'e~\eref{eq-3b}.

Certains sous-ensembles de $S$ vont jouer un r\^ole tout \`a fait particulier dans l'analyse de $\widehat s$. Ce sont les suivants : pour $d\ge 1$, soit $S_{d}$ l'ensemble des densit\'es de $S$ qui sont constantes sur une partition de $\R$ contenant au plus $d+2$ intervalles. Le ``+2" correspond au fait qu'une densit\'e constante par morceaux sur $\R$ est n\'ecessairement nulle sur les deux intervalles non born\'es de la partition ; l'entier $d$ est donc le nombre maximal d'intervalles sur lesquels $s$ est non nulle. Notons que les ensembles $(S_{d})_{d\ge 1}$ sont croissants pour l'inclusion, contiennent donc tous $S_{1}$ qui est l'ensemble des densit\'es uniformes sur un intervalle, et qu'ils sont stables par translation et changement d'\'echelle. Ils ne sont, en particulier, pas compacts. 

Le r\'esultat suivant est d\'emontr\'e dans Baraud et Birg\'e~(2016a): 
\begin{prop}
Quelle que soit la densit\'e $s\in S_{d}$,
\[
D_{n}^{\sbP}(\gP_{\!s})\le d\log_{+}^{3}\pa{n\over d}\quad \mbox{o\`u}\quad \log_{+}(u)=\max\{1,\log u\},
\]
donc, en notant $\E_{s\et}$ pour $\E_{\gP\et}$,
\[
\sup_{s\et\in S_{d}}\E_{s\et}\cro{h^{2}(s\et,\widehat s)}\le C{d\over n}\log_{+}^{3}\pa{n\over d}.
\]
\end{prop}
On a ainsi obtenu une borne de risque uniforme sur chaque ensemble $S_{d}$, 
ce qui signifie que si $s\et\in S_{d}$, $\widehat s$ converge vers $s\et$ (dans le gros espace $S$) quand $n\rightarrow+\infty$ \`a vitesse (en termes de distance de Hellinger) presque param\'etrique, c'est-\`a-dire en $1/\sqrt{n}$ au facteur logarithmique pr\`es. Notons que cette vitesse ne d\'epend que du nombre de morceaux $d$ et non de leur longueur : on estime donc aussi bien une densit\'e sur $[0,1]$ qu'une densit\'e sur $[0,10^{6}]$ ou sur $[0,10^{-6}]$, puisque l'espace $S_{d}$ est stable par changement d'\'echelle. La vitesse ne d\'epend pas davantage de la norme infinie de la densit\'e $s\et$. La vitesse est (probablement) approximativement minimax et le facteur $\log$ n\'ecessaire (mais peut-\^etre pas \`a cette puissance). En effet, il est connu que le risque minimax sur l'ensemble des densit\'es sur $[0,1]$, constantes par morceaux avec $d$ morceaux est au moins $C(d/n)\log_{+}(n/d)$ --- cf. Birg\'e et Massart~\citeyearpar{MR1653272}  --- et il est peu probable que la contrainte de monotonie permette de supprimer le $\log$.

La propri\'et\'e fondamentale de stabilit\'e exprim\'ee par \eref{eq-3b} (en prenant $\overline \gP=\gP_{\!s}$ avec $s\in S_{d}$) dit que 
\[
\E_{s\et}\cro{h^{2}(s\et,\widehat s)}\le C\cro{h^{2}(s\et,s)+ {d\over n}\log_{+}^{3}\pa{n\over d}}.
\]
Comme $d$ et $s$ dans $S_{d}$ sont arbitraires, on conclut que
\begin{equation}\label{eq4}
\E_{s\et}\cro{h^{2}(s\et,\widehat s)}\le C\inf_{d\ge 1}\cro{h^{2}(s\et,S_{d})+ {d\over n}\log_{+}^{3}\pa{n\over d}}.
\end{equation}
Donc, si $s\et$ est proche d'une densit\'e  $s\in S_{d}$ pour un certain $d$, la borne de risque pr\'ec\'edente est pratiquement la m\^eme que celle que l'on obtiendrait pour $s\et=s\in S_{d}$. 

Nous allons utiliser ce dernier r\'esultat pour aller plus loin dans l'analyse des performances de $\widehat s$. 
Pour $M\ge0$, soit $S(M)$ l'ensemble des densit\'es $s$ de $S$ \`a support sur un intervalle compact $I_{s}$ (pouvant donc d\'ependre de $s$) et telles que 
\[
\lambda(I_{s})\left[\sup_{I_{s}}s-\inf_{I_{s}}s\right]\le M.
\] 
L'ensemble $S(M)$ est encore stable par translation et changement d'\'echelle et contient toutes les lois uniformes (correspondant \`a $M=0$). On \'etablit dans Baraud et Birg\'e~(2016a) le r\'esultat d'approximation suivant. 
%
\begin{prop}\label{prop-}
Pour tout $s\in S(M)$ et $d\ge 1$
\[
h^{2}(s,S_{d})\le {M\over4d^{2}}\wedge 1.
\]
\end{prop}
Il d\'ecoule alors de~\eref{eq4} que, si $s\et\in S(M)$,
\begin{align*}
\E_{s\et}\cro{h^{2}(s\et,\widehat s)}&\le C\inf_{d\ge 1}\cro{\left({M\over4d^{2}}\wedge 1\right)+{d\over n}\log_{+}^{3}\pa{n\over d}}\\&\le C'\max\left\{{M^{1/3}\over n^{2/3}}\log^{2}n,{\log_{+}^{3}(n)\over n}\right\},
\end{align*}
ce qui conduit \`a des bornes de risque de $\widehat s$ uniformes sur les sous-ensembles $S(M)\subset S$. Si $s\et$ n'est pas dans $S(M)$ mais est proche d'un \'el\'ement $s$ de $S(M)$ alors le risque en $s\et$ est major\'e par
\[
\E_{s\et}\cro{h^{2}(s\et,\widehat s)}\le C\cro{h^{2}(s\et,s)+ \max\left\{{M^{1/3}\over n^{2/3}}\log^{2}n,{\log_{+}^{3}(n)\over n}\right\}}.
\]
Comme $M$ et $s$ dans $S(M)$ et sont arbitraires, on en d\'eduit que
\begin{equation}
\E_{s\et}\cro{h^{2}(s\et,\widehat s)}\le C\inf_{M\ge 0}\cro{h^{2}(s\et,S(M))+ \max\left\{{M^{1/3}\over n^{2/3}}\log^{2}n,{\log_{+}^{3}(n)\over n}\right\}}.
\label{eq7}
\end{equation}
On peut poursuivre l'analyse de notre borne de risque pour les $\rho$-estimateurs si l'on a de l'information sur la mani\`ere dont les espaces $S(M)$ approximent la densit\'e inconnue $s\et$. Supposons, par exemple, que $s\et(x)=\theta \exp\cro{-\theta x}\1_{\R_{+}}(x)$ pour un certain $\theta>0$, donc $s\et\in S$. Pour tout $T>0$, on peut l'approximer par la densit\'e
\[
\overline s_{T}(x)={\theta e^{-\theta x}\over 1-e^{-\theta T}}\1_{[0,T]}.
\]
Comme 
\[
T\left[\sup_{[0,T]}\overline{s}_{T}-\inf_{[0,T]}\overline{s}_{T}\right]=T\cro{{\theta\over 1-e^{-\theta T}}-{\theta e^{-\theta T}\over 1-e^{-\theta T}}}=\theta T,
\]
$\overline s_{T}\in S(\theta T)$. De plus,
\[
\rho(P_{s\et},P_{\overline s_{T}})=\int_{0}^{T}{\theta e^{-\theta x}\over \sqrt{1-e^{-\theta T}}}dx= \sqrt{1-e^{-\theta T}},
\]
donc 
\[
h^{2}(s\et,\overline{s}_{T})=1-\sqrt{1-e^{-\theta T}}={e^{-\theta T}\over 1+\sqrt{1-e^{-\theta T}}}\le e^{-\theta T}.
\]
Il d\'ecoule alors de \eref{eq7} que, puisque $T$ est arbitraire, donc $M=\theta T$ \'egalement,
\[
\E_{s\et}\cro{h^{2}(s\et,\widehat s)}\le C\inf_{M\ge 0}\cro{e^{-M}+\max\left\{{M^{1/3}\over n^{2/3}}\log^{2}n,{\log_{+}^{3}(n)\over n}\right\}}.
\]
En choisissant $M=(2/3)\log n$, on obtient finalement une borne de risque 
\[
\E_{s\et}\cro{h^{2}(s\et,\widehat s)}\le C{\log^{7/3} n\over n^{2/3}}
\]
qui est ind\'ependante de $\theta>0$, ce qui montre que $\widehat s$ converge \`a vitesse ${n^{-1/3}\log^{7/6} n}$ uniform\'ement sur l'ensemble de toutes les lois exponentielles bien que l'estimateur n'ait pas \'et\'e particuli\`erement con\c cu pour estimer de telles lois. Le r\'esultat demeure valable pour les lois exponentielles translat\'ees de densit\'es
$\theta \exp\cro{-\theta(x-a)}\1_{[a,+\infty)}(x)$ avec $a\in\R$.

La d\'emarche pr\'ec\'edente n'est pas sp\'ecifique \`a l'ensemble $S$ et l'on peut consid\'erer d'autres exemples de familles de densit\'es d\'efinies par des contraintes de forme et pour lesquelles les choses se passent \`a peu pr\`es de la m\^eme mani\`ere. On peut en particulier \'etudier --- cf. Baraud et Birg\'e~(2016a) --- :
\begin{lista}
\item l'ensemble des densit\'es monotones sur chaque \'el\'ement d'une partition de $\R$ en au plus $k$ intervalles;
\item  l'ensemble des densit\'es dont la racine carr\'ee est convexe ou concave sur chaque \'el\'ement d'une partition de $\R$ en au plus $k$ intervalles;
\item l'ensemble des densit\'es qui sont log-concaves sur $\R$, c'est-\`a-dire de la forme $e^{g}\1_{I}$ o\`u $I$ est un intervalle de $\R$ et $g$ une fonction concave. Par exemple, les lois gaussiennes, expontielles, de Laplace ou uniformes sont log-concaves. Dans ce cas ce sont les densit\'es de la forme $e^{g}\1_{I}$ avec $g$ lin\'eaire par morceaux qui vont jouer le r\^ole des fonctions constantes par morceaux de l'exemple pr\'ec\'edent. En particulier, 
on peut montrer que le $\rho$-estimateur d\'efini sur ce gros espace de densit\'es va converger \`a vitesse param\'etrique (\`a des facteurs logarithmiques pr\`es) vers toutes les lois uniformes, exponentielles, ou de Laplace (pour lesquelles $g$ est effectivement lin\'eaire par morceaux).
\end{lista}
%

\subsection{Le mod\`ele de r\'egression \`a plan al\'eatoire\label{S10b}}

\subsubsection{Estimation sur un mod\`ele\label{S10b1}}
Nous supposons ici que les observations $X_{i}$, $i\in\{1,\ldots,n\}$ sont ind\'ependantes et de la forme $(W_{i},Y_{i})\in \sW\times \R$ et que notre mod\`ele statistique est d\'efini sur la base des hypoth\`eses suivantes : 
\begin{listi}
\item les $W_{i}$ sont i.i.d.\ de loi $P_{W}$ inconnue;
\item il existe une fonction $f\et$ dans un espace vectoriel $\sF$ de dimension $d$  telle que 
\begin{equation}
Y_{i}=f\et(W_{i})+\eps_{i}\quad \mbox{pour tout $i\in\{1,\ldots,n\}$},
\label{eq2}
\end{equation}
et les v.a.r.\ $\eps_{i}$ sont i.i.d.\ de densit\'e unimodale $q$ par rapport \`a la mesure de Lebesgue $\lambda$ sur $\R$.
\end{listi}
Ceci signifie que le mod\`ele que l'on va utiliser pour estimer la loi $\gP\et$ de $\biX$ est de la forme $\sbP=\{\gP_{\!f}=P_{f}^{\otimes n},\ f\in \sF\}$ avec, pour tout $f\in\sF$,
\[
P_{f}=q_{f}\cdot (P_{W}\otimes \lambda)\qquad \mbox{et} \qquad q_{f}(w,y)=q(y-f(w)).
\] 
Comme le crit\`ere permettant de calculer le $\rho$-estimateur $\widehat \gP=P_{\widehat f}^{\otimes n}$ ne d\'epend que des rapports $q_{f}/q_{f'}$ pour $f,f'\in \sF$, il n'est pas n\'ecessaire de connaitre la loi $P_{W}$ des $W_{i}$ pour le calculer. 
 
Sous l'hypoth\`ese que $q$ est unimodale et que $\sF$ un espace vectoriel de dimension $d\le n$, on peut montrer que 
\[
D_{n}^{\sbP}(\overline \gP)\le Cd\log(en/d)\quad \mbox{pour tout }\overline \gP\in \sbP.
\]
On peut m\^eme remplacer l'hypoth\`ese ``unimodale" par ``$k$-modale" (qui a au plus $k$ modes) et $\sF$ par une classe de fonctions VC-subgraph d'indice $d\ge 1$ et obtenir la m\^eme in\'egalit\'e avec une constante $\overline{C}$ d\'ependant alors de $k$.

Dans cette situation, il d\'ecoule de~\eref{eq-3b} que
\begin{equation}
\E_{\gP\et}\cro{H^{2}\left(\gP\et,\gP_{\!\widehat f}\right)}\le C\left[\inf_{f\in \sF}H^{2}(\gP\et,\gP_{\!f})+d\log\pa{en\over d}\right].
\label{eq-rand}
\end{equation}
Pour analyser cette in\'egalit\'e, notamment ce qu'elle dit de l'estimation d'une fonction de r\'egression, et la comparer aux r\'esultats habituels, il convient de supposer que le mod\`ele de r\'egression \eref{eq2} est exact \`a ceci pr\`es que $f\et$ n'appartient pas n\'ecessairement \`a $\sF$ et que les $\eps_{i}$ n'ont pas n\'ecessairement la densit\'e $q$ mais une densit\'e $p$, ce qui revient \`a dire que, pour tout $i\in\{1,\ldots,n\}$,
\[
P_{i}\et=p_{f\et}\cdot (P_{W}\otimes \lambda) \qquad \mbox{avec} \qquad p_{f\et}(w,y)=p(y-f\et(w)).
\]
Pour relier ais\'ement la distance $H$ \`a une distance entre les fonctions de r\'egression nous ferons maintenant l'hypoth\`ese que $f\et$ ainsi que tous les \'el\'ements de $\sF$ sont born\'es par une constante $B$ et que la densit\'e $q$ est telle que
\[
|u|^{-r}\left[{1\over 2}\int_{\R}\pa{\sqrt{q(y)}-\sqrt{q(y-u)}}^{2}dy\right]
\fleche_{u\rightarrow0}A\in(0,+\infty)\quad \mbox{pour un }r\in (0,2],
\]
c'est-\`a-dire que le carr\'e de la distance de Hellinger entre les densit\'es $q$ et $q(\cdot-u)$ est \'equivalent \`a $A|u|^{r}$ lorsque $u\rightarrow0$.
On sait que si le mod\`ele de translation associ\'e \`a $q$ est r\'egulier (densit\'es gaussienne, de Cauchy, de Laplace, etc.) alors $r=2$ et que $r=1$ si $q$ est la densit\'e uniforme sur un intervalle. Dans ces conditions, on peut montrer que, pour tout $f\in\sF$,
\begin{equation}
\underline{c}\left[\norm{f\et-f}_{r}^{r}+h^{2}(p,q)\right]\le \frac{H^{2}\left(\gP\et,\gP_{\!f}\right)}{n}\le \overline{c}\left[\norm{f\et-f}_{r}^{r}+h^{2}(p,q)\right]
\label{eq3}
\end{equation}
avec des quantit\'es $\underline{c}$ et $\overline{c}$ d\'ependant de nos hypoth\`eses et
\[
\norm{f\et-f}_{r}^{r}=\int_{\sW}\ab{f\et(w)-f(w)}^{r}dP_{W}(w).
\]
Comme les paires $(f,q)$ et $(f\et,p)$ jouent des r\^oles sym\'etriques pour le calcul de 
$H^{2}\left(\gP\et,\gP_{\!f}\right)$, (\ref{eq3}) demeure vraie si, pour un $r\in (0,2]$,
\begin{equation}
|u|^{-r}\left[{1\over 2}\int_{\R}\pa{\sqrt{p(y)}-\sqrt{p(y-u)}}^{2}dy\right]
\fleche_{u\rightarrow0}A\in(0,+\infty).
\label{eq-ur}
\end{equation}
%
\begin{thm}\label{thm-regrand}
Si la relation (\ref{eq3}) est satisfaite, alors
\begin{equation}
\E_{\gP\et}\left[\norm{f\et-\widehat{f}}_{r}^{r}\right]\le \overline{C}\left[\inf_{f\in\sbF}\norm{f\et-f}_{r}^{r}+h^{2}(p,q)+{d\over n}\log\pa{en\over d}\right].
\label{eq-4}
\end{equation}
\end{thm}
Ce r\'esultat appelle quelques remarques:
\begin{lista}
\item il ne suppose rien sur la loi $P_{W}$ des $W_{i}$;
\item il ne suppose rien non plus sur l'int\'egrabilit\'e des erreurs et l'on peut parfaitement prendre par exemple la loi de Cauchy pour mod\'eliser la loi des $\varepsilon_{i}$, de m\^eme que la vraie densit\'e $p$ des erreurs peut \^etre la loi de Cauchy;
\item les termes $\inf_{f\in\sF}\norm{f\et-f}_{r}^{r}$ et $h^{2}(p,q)$ montrent que la borne reste stable \`a une possible erreur de sp\'ecification du mod\`ele.
\end{lista}
Pour analyser la borne de risque fournie par le Th\'eor\`eme~\ref{thm-regrand}, nous allons supposer que $p=q$, $f\et\in\sF$ (c'est-\`a-dire que notre mod\`ele est exact) et que $\sF$ est une partie d'un sous-espace vectoriel $\overline{\sF}$ de dimension $d$ engendr\'e par les fonctions $\varphi_{1},\ldots,\varphi_{d}$ de sorte que
\[
f\et=\sum_{j=1}^{d}\beta_{j}\et\varphi_{j}\qquad\mbox{et}\qquad
\widehat{f}=\sum_{j=1}^{d}\widehat{\beta}_{j}\varphi_{j}.
\]
Dans ce cas, en vertu de (\ref{eq-4}), 
\[
\E_{\gP\et}\left[\norm{f\et-\widehat{f}}_{r}^{r}\right]\le \overline{C}{d\over n}\log\pa{en\over d}.
\]
Donc, par l'in\'egalit\'e de Markov, avec une probabilit\'e proche de 1,
\[
\norm{f\et-\widehat{f}}_{r}\le \overline{C}\left[{d\over n}\log\pa{en\over d}\right]^{1/r}
\]
et, comme sur l'espace vectoriel $\overline{\sF}$ toutes les normes sont \'equivalentes,
cette relation implique que, avec une probabilit\'e proche de 1,
\[
\max_{1\le j\le d}\left|\widehat{\beta}_{j}-\beta_{j}\et\right|\le \overline{C}\left[\frac{\log n}{n}\right]^{1/r},
\]
o\`u la quantit\'e $\overline{C}$ d\'epend aussi de $d$. Si $r=2$ on retrouve une vitesse de convergence param\'etrique classique en $1/\sqrt{n}$ (au facteur logarithmique pr\`es) mais la vitesse de convergence est plus rapide si $r<2$. Par exemple, dans le cas d'erreurs de loi uniforme pour lequel $r=1$, le $\rho$-estimateur converge \`a vitesse $1/n$ (au facteur $\log$ pr\`es) quand l'estimateur des moindres carr\'es ordinaire converge lui \`a vitesse $1/\sqrt{n}$. Si la densit\'e des erreurs n'est pas born\'ee, $r<1$ et la convergence est encore plus rapide. 

Nous insisterons ici sur le fait que les hypoth\`eses que nous avons utilis\'ees pour obtenir
la borne g\'en\'erale (\ref{eq-rand}) sont extr\^emement faibles : aucune pour la loi $P_{W}$ ni la structure de l'espace vectoriel $\sbF$ et une simple hypoth\`ese d'unimodalit\'e sur la densit\'e des erreurs $\varepsilon_{i}$ (que l'on peut d'ailleurs remplacer par une borne sur le nombre de modes de cette densit\'e). Ceci contraste avec les hypoth\`eses que l'on rencontre d'ordinaire pour traiter ce probl\`eme, en particulier pour contr\^oler le risque des estimateurs des moindres carr\'es. En contrepartie, on peut regretter la pr\'esence du facteur logarithmique dans nos bornes de risque puisque, dans certaines situations, on obtient, mais au prix d'hypoth\`eses nettement plus fortes, des vitesses analogues sans le facteur logarithmique. En contrepartie, notre estimateur est robuste et ses performances ne seront que l\'eg\`erement affect\'ees par un petit nombre (petit devant $n$) d'observations $(W_{i},Y_{i})$ aberrantes. A contrario, une seule paire $(W_{i},Y_{i})$ suffisamment \og exotique \fg{} peut compl\`etement d\'er\'egler le comportement des estimateurs des moindres carr\'es.

\subsubsection{Plusieurs mod\`eles}
Il est \'evidemment assez restrictif de ne travailler qu'avec une seule densit\'e $q$ et un seul mod\`ele $\sF$ mais, comme indiqu\'e \`a la section~\ref{S7}, nous pouvons en fait utiliser plusieurs mod\`eles simultan\'ement de mani\`ere \`a faire varier la densit\'e $q$ ainsi que l'espace $\sF$. Nous pouvons consid\'erer des $\rho$-estimateurs construits \`a partir de familles d\'enombrables de densit\'es $q$ pour mod\'eliser $p$ et plusieurs familles d'espaces fonctionnels $\sF$ pour mod\'eliser $f\et$. Nous pouvons, par exemple, utiliser une famille $\{\sF_{d},d\in\N\setminus\{0\}\}$ o\`u $\sF_{d}$ est un sous-ensemble d'un espace vectoriel $\overline{\sF}_{\!d}$ de dimension $d$ ainsi qu'une famille d\'enombrable $\{q_{k},k\in\cK\subset\N\}$ de densit\'es, ce qui fournit une famille de mod\`eles de la forme $(\sF_{d},q_{k})$ index\'es par $(d,k)$, un tel mod\`ele correspondant \`a l'hypoth\`ese (a priori inexacte) que
\[
Y_{i}=f\et(W_{i})+\eps_{i}\quad \mbox{pour tout $i\in\{1,\ldots,n\}$},\quad f\et\in\sF_{d}\quad\mbox{et}\quad
\varepsilon_{i}\sim q_{k}\cdot\lambda.
\]
Dans ce cadre, un $\rho$-estimateur aura la forme $\gP_{\!(\widehat{f},\widehat{q})}$ avec $\widehat{f}\in\bigcup_{d\ge1}\sF_{d}$ et $\widehat{q}=q_{\widehat{k}}$, $\widehat{k}\in\N$.
Comme nous l'avons vu, pour un mod\`ele $(\sF_{d},q_{k})$,
\[
D_{n}^{\sbP}(\overline \gP)\le C_{0}d\log(en/d)\quad \mbox{pour tout }\overline \gP\in(\sF_{d},q_{k}).
\]
Soit alors une suite $(a_{k})_{k\in\cK}$ de nombres positifs tels que $\sum_{k\in\cK}e^{-a_{k}}=e-1$.
Si nous fixons $\Delta(d,k)=d+a_{k}$, nous trouvons que
\begin{align*}
\sum_{d=1}^{+\infty}\sum_{k\in\cK}e^{-\Delta(d,k)}&=\sum_{d=1}^{+\infty}\sum_{k\in\cK}e^{-(d+a_{k})}=\left(\sum_{d=1}^{+\infty}e^{-d}\right)\left(\sum_{k\in\cK}e^{-a_{k}}\right)
\\&={e^{-1}\over 1-e^{-1}}(e-1)=1.
\end{align*}
Nous pouvons donc fixer la p\'enalit\'e de la mani\`ere suivante :
\[
\pen(\gQ)=\kappa\left[C_{0}d\log(en/d)+d+a_{k}\right]\quad\mbox{pour tout }\gQ\in(\sF_{d},q_{k}).
\]
Un $\rho$-estimateur $(\widehat{f},\widehat{q})$ construit sur une telle famille de mod\`eles aura alors un risque born\'e par
\[
\E_{\gP\et}\cro{H^{2}\left(\gP\et,\gP_{\!(\widehat{f},\widehat{q})}\right)}\le C\inf_{(d,k)}\left[\inf_{f\in \sF_{d}}H^{2}(\gP\et,\gP_{\!(f,q_{k})})+d\log\pa{en\over d}+a_{k}\right].
\]
Notons que, si $f\et$ et tous les \'el\'ements de $\bigcup_{d\ge1}\sF_{d}$ sont uniform\'ements born\'es par $B$ et si la densit\'e $p$ v\'erifie la condition (\ref{eq-ur}), les in\'egalit\'es (\ref{eq3}) seront satisfaites pour tout $f\in\bigcup_{d\ge1}\sF_{d}$, ce qui conduira \`a une borne de risque de la forme
\begin{align*}
\lefteqn{\E_{\gP\et}\left[\norm{f\et-\widehat{f}}_{r}^{r}+h^{2}(p,\widehat{q})\right]}\hspace{10mm}\\&\le \overline{C}\left[\inf_{d}\left(\inf_{f\in\sbF_{d}}\norm{f\et-f}_{r}^{r}+{d\over n}\log\pa{en\over d}\right)+\inf_{k}\left[h^{2}(p,q_{k})+a_{k}\right]\right].
\end{align*}

\section{Conclusion}
L'int\'er\^et principal des $\rho$-estimateurs r\'eside dans leur robustesse, comme le montre la borne \eref{eq-1b}, ce qui permet de travailler syst\'ematiquement avec des mod\`eles approch\'es et de se prot\'eger contre la pr\'esence \'eventuelle d'observations atypiques. Cette robustesse permet \'egalement de remplacer un mod\`ele complexe par un ou plusieurs mod\`eles plus simples de mani\`ere \`a optimiser le compromis entre l'erreur d'approximation et l'erreur d'estimation li\'ee \`a la complexit\'e du mod\`ele sur lequel l'estimateur est construitet d'atteindre ainsi la vitesse optimale d'estimation (\'eventuellement \`a un facteur logarithmique pr\`es). De plus le $\rho$-estimateur s'appuie sur des mod\`eles de probabilit\'es et ses performances ne d\'ependent ni de la mesure dominante ni du choix des densit\'es. 

Comme nous l'avons vu dans la section~\ref{EMV}, toutes ces qualit\'es ne sont pas partag\'ees par l'estimateur du maximum de vraisemblance ou l'estimateur des moindres carr\'es en r\'egression. En revanche le $\rho$-estimateur, qui co\"{\i}ncide (asymptotiquement ou non) avec l'estimateur du maximum de vraisemblance dans un certain nombre de situations, b\'en\'eficie alors de ses propri\'et\'es d'optimalit\'e.

Dans la section~\ref{ex-sect2.5}, nous avons mis en \'evidence un probl\`eme li\'e \`a l'utilisation du maximum de vraisemblance sur un mod\`ele approch\'e, nous allons voir \`a pr\'esent, en guise de conclusion, ce qui se passe si l'on remplace le maximum de vraisemblance par un $\rho$-estimateur construit sur ce m\^eme mod\`ele approch\'e.

\paragraph{Retour sur l'exemple de la section~\ref{ex-sect2.5}.}
Nous avons vu que la m\'ethode du maximum de vraisemblance sur le mod\`ele  approchant $\sP_{0}$ ne fournissait, avec une probabilit\'e proche de un, aucune estimation de $\theta$. Une question naturelle est de savoir ce qui se passe si nous utilisons un $\rho$-estimateur et quelles sont alors ses performances. Si l'on note $p_{\theta}$ la densit\'e $ \1_{[\theta,\theta+1]}(x)$ de la loi uniforme sur $[\theta,\theta+1]$, le $\rho$-estimateur $\widehat \theta_{n}$ de $\theta$ minimise sur $\R$ l'application
\begin{equation}\label{eq-crit}
\theta\mapsto \sup_{\theta'\in\R}\sum_{i=1}^{n}\psi\pa{\sqrt{p_{\theta'}(X_{i})\over p_{\theta}(X_{i})}}
\end{equation}
avec  $\psi(1/0)=\psi(+\infty)=1=-\psi(0)$ et $\psi(1)=0$. Or, pour tout $\theta,\theta'\in\R$ et $x\in\R$, 
\[
\psi\pa{\sqrt{p_{\theta'}(x)\over p_{\theta}(x)}}=
\begin{cases}
1 &\mbox{si}\quad  x\in [\theta',\theta'+1]\setminus [\theta,\theta+1],\\
0 &\mbox{si}\quad  x\in [\theta',\theta'+1]\cap [\theta,\theta+1],\\
-1 &\mbox{si}\quad  x\in [\theta,\theta+1]\setminus [\theta',\theta'+1],
\end{cases}
\]
donc $\psi\pa{\sqrt{p_{\theta'}(x)/p_{\theta}(x)}}=\1_{[\theta',\theta'+1]}(x)-\1_{[\theta,\theta+1]}(x)$ et 
\[
\sum_{i=1}^{n}\psi\pa{\sqrt{p_{\theta'}(X_{i})\over p_{\theta}(X_{i})}}=\sum_{i=1}^{n}\1_{[\theta',\theta'+1]}(X_{i})-\sum_{i=1}^{n}\1_{[\theta,\theta+1]}(X_{i}).\\
\]
Alors 
\[
\sup_{\theta'\in\R}\sum_{i=1}^{n}\psi\pa{\sqrt{p_{\theta'}(X_{i})\over p_{\theta}(X_{i})}} =\cro{ \sup_{\theta'\in\R}\sum_{i=1}^{n}\1_{[\theta',\theta'+1]}(X_{i})}-\sum_{i=1}^{n}\1_{[\theta,\theta+1]}(X_{i})
\]
et un $\rho$-estimateur $\widehat \theta_{n}$ est n'importe quel \'el\'ement maximisant l'application
\[
\theta\mapsto \sum_{i=1}^{n}\1_{[\theta,\theta+1]}(X_{i}).
\] 
En d'autres termes, calculer un $\rho$-estimateur revient \`a rechercher un param\`etre $\theta$ maximisant le nombre de donn\'ees $X_{i}$ tombant dans l'intervalle $[\theta,\theta+1]$ quand l'e.m.v.\ recherche un param\`etre $\theta$ pour lequel {\em toutes les donn\'ees} $X_{1},\ldots,X_{n}$ appartiennent \`a $[\theta,\theta+1]$. Mais un tel param\`etre n'existe pas d\`es que l'on observe une valeur aberrante.

\'Etudions \`a pr\'esent les performances de $\widehat \theta_{n}$. Soit $\theta_{0}\in \R$ et $\alpha\in [0,\eps]$ avec $\eps<1/2$ les vraies valeurs des param\`etres de la loi de l'observation $\bsX=(X_{1},\ldots,X_{n})$. Notre objectif est d'\'evaluer la probabilit\'e 
\[
\P_{\alpha,\theta_{0}}\cro{|\theta_{0}-\widehat \theta_{n}|>c/n}\quad \mbox{pour $0<c<n$.}
\]
Notons que sous $P_{\alpha,\theta_{0}}$, la loi jointe des variables $X_{i}-\theta_{0}$ est ind\'ependante de $\theta_{0}$ et la valeur d'un $\rho$-estimateur  $\widehat \theta_{n}=\widehat \theta_{n}(X_{1},\ldots,X_{n})$ construit \`a partir des observations $X_{1},\ldots,X_{n}$ s'\'ecrit sous la forme $\theta_{0}+\widetilde \theta_{n}$ o\`u $\widetilde \theta_{n}$ est un $\rho$-estimateur construit \`a partir des observations $X_{1}-\theta_{0},\ldots,X_{n}-\theta_{0}$ de sorte qu'il nous suffit d'\'etudier le cas $\theta_{0}=0$. 

Sous $P_{\alpha,0}$,  le nombre $N_{n}$ de $X_{i}$ appartenant \`a l'intervalle $[0,1]$ suit une loi binomiale de param\`etres $n$ et $1-\alpha$ et deux situations peuvent alors se produire: soit $N_{n}\le n/2$, ce qui arrive avec une probabilit\'e 
\[
\P_{\alpha,0}\cro{N_{n}\le n/2}=\P_{\alpha,0}\cro{N_{n}-n(1-\alpha)\le-(n/2)(1-2\alpha)}\le \exp\cro{-n(1-2\alpha)^{2}/2},
\] 
d'apr\`es l'in\'egalit\'e de Hoeffding, soit $N_{n}>n/2$ auquel cas un $\rho$-estimateur $\widehat \theta_{n}$ est un point arbitraire de $[X_{(n)}-1,X_{(1)}]$ o\`u $X_{(1)}=\min_{i=1,\ldots,n} X_{i}$ et $X_{(n)}$ est le plus grand des $X_{i}$ appartenant \`a l'intervalle $[0,1]$. Sur l'\'ev\`enement $\{N_{n}>n/2\}$, $X_{(n)}=\max_{i=1,\ldots,n}X_{i}\1_{X_{i}\in [0,1]}$.
Par des calculs classiques, 
\[
\P_{\alpha,0}\cro{X_{(1)}>c/n}=\left(\P_{\alpha,0}\cro{X_{1}>c/n}\st\right)^{n}=\cro{1-{c(1-\alpha)\over n}}^{n}\le e^{-(1-\alpha)c}
\]
et 
\begin{align*}
\P_{\alpha,0}\cro{X_{(n)}-1<-c/n,\,N_{n}>n/2}&\le \left(\P_{\alpha,0}\cro{X_{1}\1_{[0,1]}(X_{1})<1-c/n}\right)^{n}\\
&\le\left(\P_{\alpha,0}\cro{\{0\le X_{1}<1-c/n\}\cup\{X_{1}>1\}}\right)^{n}\\
&=\cro{1-{(1-\alpha)c\over n}}^{n}\le e^{-(1-\alpha)c}.
\end{align*}
Ainsi, pour tout $\rho$-estimateur $\widehat \theta_{n}$, 
\begin{align*}
\P_{\alpha,0}\cro{|\widehat \theta_{n}|>c/n}
\le\:& \P_{\alpha,0}\cro{N_{n}\le n/2}+\P_{\alpha,0}\cro{N_{n}>n/2,|\widehat \theta_{n}|>c/n}\\
\le\:&\P_{\alpha,0}\cro{N_{n}\le n/2}+\P_{\alpha,0}\cro{X_{(1)}>c/n}\\&+\P_{\alpha,0}\cro{N_{n}>n/2,X_{(n)}-1<-c/n}\\
\le\:& e^{-n(1-2\alpha)^{2}/2}+2e^{-(1-\alpha)c}.
\end{align*}

Ce r\'esultat montre que pour toute loi $P_{\alpha,\theta_{0}}\in \sP_{\eps}$, le $\rho$-estimateur de $\theta$ construit sur le sous-mod\`ele $\sP_{0}$ converge \`a la vitesse $1/n$  d\`es que $\eps<1/2$ alors m\^eme que, si $\alpha$ n'est pas nul, la distance de Hellinger entre la vraie loi $P_{\alpha,\theta_{0}}$ et le mod\`ele $\sP_{0}$ est strictement positive et ind\'ependante de $n$. En plus d'illustrer la robustesse du $\rho$-estimateur, cet exemple montre qu'une in\'egalit\'e telle que \eref{eq-RisqueMax} sur le risque du $\rho$-estimateur en distance de Hellinger peut s'av\'erer pessimiste dans des cadres param\'etriques pour lesquels notre objectif n'est pas  d'estimer la loi de l'observation mais simplement son param\`etre pour la perte euclidienne usuelle.

\bibliographystyle{apalike}

\end{document}